\documentclass[reqno]{amsart}

\usepackage{float}
\usepackage[utf8]{inputenc}
\usepackage{comment}
\usepackage{capt-of}
\usepackage{fullpage}
\usepackage[usenames,dvipsnames,svgnames,table]{xcolor}
\usepackage[colorlinks=true, pdfstartview=FitV, linkcolor=blue, citecolor=blue, urlcolor=darkblue]{hyperref}

\usepackage[maxbibnames=15]{biblatex}
\usepackage[margin=1in]{geometry}                
\usepackage{graphicx}
\usepackage{mathrsfs}
\usepackage{amssymb}
\usepackage{amsmath}
\usepackage{amsfonts}
\usepackage{mathrsfs}
\usepackage{epstopdf}
\usepackage{lscape}
\usepackage{array}
\usepackage[utf8]{inputenc}
\usepackage{tikz,caption}
\usetikzlibrary{arrows.meta}
\usepackage[capitalise]{cleveref}

\usepackage{enumitem}
\setlist[itemize]{leftmargin=2em}
\setlist[enumerate]{leftmargin=2em}

\usepackage{booktabs}
\usepackage{multirow}
\usepackage{mathtools}
\usepackage[linesnumbered,ruled]{algorithm2e}
\usepackage{tikz}
\usepackage{mathtools}
\usepackage{soul}
\usepackage{upgreek}

\newtheorem{theorem}{Theorem}[section]
\newtheorem{lemma}[theorem]{Lemma}
\newtheorem{corollary}[theorem]{Corollary}
\newtheorem{conjecture}[theorem]{Conjecture}
\newtheorem{proposition}[theorem]{Proposition}

\theoremstyle{definition}
\newtheorem{definition}[theorem]{Definition}
\newtheorem{notation}[theorem]{Notation}
\newtheorem{example}[theorem]{Example}

\newtheorem{remark}[theorem]{Remark}

\newcolumntype{C}[1]{>{\centering\arraybackslash}m{#1}}

\newcommand{\Q}{\mathbb{Q}}

\newcommand{\St}{\mathsf{St}}

\newcommand{\mf}[1]{\mathfrak{#1}}
\newcommand{\mc}[1]{\mathcal{#1}}
\renewcommand{\ul}[1]{\underline{#1}}

\newcommand{\la}{\lambda}

\newcommand{\stf}{\mathfrak{st}}

\newcommand{\sort}{\mathrm{sort}}
\newcommand{\lc}{\lambda_\mathrm{LC}}
\newcommand{\lead}{{\lambda_{\mathrm{lead}}}}

\DeclarePairedDelimiter{\abs}{\lvert}{\rvert}

\definecolor{darkblue}{rgb}{0.0,0,0.7} % darkblue color
\definecolor{darkred}{rgb}{0.7,0,0} % darkred color
\definecolor{darkgreen}{rgb}{0, .6, 0} % darkgreen color

\newcommand{\defncolor}{\color{darkred}}
\newcommand{\defn}[1]{{\defncolor\emph{#1}}} % emphasis of a definition

\usepackage[colorinlistoftodos]{todonotes}

\title{The Chromatic Symmetric Function for Unicyclic Graphs}
\author[A.~Bingham]{Aram Bingham}\address[A.~Bingham]{Departamento de Matem\'aticas, Universidad de Chile, 3425 Las Palmeras, Ñuñoa, Santiago, CL}
\email{\textcolor{blue}{\href{mailto:aram@matmor.unam.mx}{aram@matmor.unam.mx}}}

\author[L.~Johnston]{Lisa Johnston}
\address[L.~Johnston]{Department of Mathematics, University of California, One Shields Avenue, Davis, CA, USA}
\email{\textcolor{blue}{\href{mailto:lisjohnston@ucdavis.edu}{lisjohnston@ucdavis.edu}}}

\author[C.~Lawson]{Colin Lawson}
\address[C.~Lawson]{Department of Mathematics \& Statistics, Stephen F. Austin State University, Nacogdoches, TX}
\email{\textcolor{blue}{\href{mailto:Colin.Lawson@sfasu.edu}{Colin.Lawson@sfasu.edu}}}

\author[R.~Orellana]{Rosa Orellana}
\address[R.~Orellana]{Mathematics Department, Dartmouth College,
Hanover, NH, USA}
\email{Rosa.C.Orellana@dartmouth.edu}
\urladdr{\href{https://math.dartmouth.edu/~orellana/}{https://math.dartmouth.edu/~orellana/}}

\author[J.~Pan]{Jianping Pan}
\address[J.~Pan]{School of Mathematical and Statistical Sciences, Arizona State University, Tempe, AZ, USA}
\email{\textcolor{blue}{\href{mailto:jianping.pan@asu.edu}{jianping.pan@asu.edu}}}

\author[C.~Sato]{Chelsea Sato}
\address[C.~Sato]{Department of Mathematics, Willamette University, Salem, OR, USA }
\email{\textcolor{blue}{\href{mailto:cmsato@willamette.edu}{cmsato@willamette.edu}}}

\date{\today}

\begin{document}

\begin{abstract}
    Motivated by the question of which structural properties of a graph can be recovered from the chromatic symmetric function (CSF), we study the CSF of connected unicyclic graphs. 
     While it is known that there can be non-isomorphic unicyclic graphs with the same CSF, we find experimentally that such examples are rare 
     for graphs with up to 17 vertices.
    In fact, in many cases we can recover data such as the number of leaves, number of internal edges, cycle size, and number of attached non-trivial trees, by extending known results for trees to unicyclic graphs. 
    These results are obtained by analyzing the CSF of a connected unicyclic graph in the \emph{star-basis} using the deletion-near-contraction (DNC) relation developed by Aliste-Prieto, Orellana and Zamora, and computing the ``leading" partition, its coefficient, as well as coefficients indexed by hook partitions. 
    We also give explicit formulas for star-expansions of several classes of graphs, developing methods for extracting coefficients using structural properties of the graph.

\end{abstract}

\maketitle
\section{Introduction}
The chromatic symmetric function 
$\mathbf{X}_G$ 
of a graph $G$ was introduced in \cite{stanley_1995} by Stanley as a symmetric function generalization of 
the chromatic polynomial of a graph. 
Specifically, for a graph $G$ with vertex set $V=\{v_1, \dots, v_n\}$,  $\mathbf{X}_G$ is defined, for commuting variables $x_1, x_2, \dots$, by
\[
\mathbf{X}_G =
\sum_\kappa x_{\kappa(v_1)}x_{\kappa(v_2)}\dots x_{\kappa(v_n)},
\]
where $\kappa: V \to \mathbb{N}$ ranges over all proper colorings of the vertices of $G$.
Setting the first $k$ variables to $1$ and all other variables to zero indeed recovers the value $\chi_G(k)$ of the chromatic polynomial of $G$.
A significant body of research on $\mathbf{X}_G$ has focused on Stanley's ``tree isomorphism problem", which asks whether the chromatic symmetric function distinguishes non-isomorphic trees. While this question remains open in general, the chromatic symmetric function has been shown to distinguish various subclasses of trees (see, e.g., \cite{ALISTEPRIETO,qCats,HurynChmutov20,LoeblSereni, martin2007distinguishingtreeschromaticsymmetric,WANG}), as well as for all trees up to 29 vertices \cite{HeilJi18}.

Many key structural properties of a graph can be recovered from $\mathbf{X}_G$ by examining its coefficients when expanded with respect to various bases of symmetric functions. For example, expanding $\mathbf{X}_G$ in the power sum symmetric function basis reveals structural properties such as the number of vertices, edges, girth and number of connected components \cite{martin2007distinguishingtreeschromaticsymmetric}, the degree sequence and number of leaves for trees \cite{martin2007distinguishingtreeschromaticsymmetric}, and the number of triangles in $G$ \cite{ORELLANA20141}. In particular, $\mathbf{X}_G$ can tell us whether $G$ is a connected unicyclic graph~\cite{martin2007distinguishingtreeschromaticsymmetric}. Recently, the chromatic symmetric function has been studied~\cite{AMOZ23,gonzalez2024chromatic} in terms of the \textit{star-basis}, $\{\mathfrak{st}_{\lambda} \mid \lambda \vdash n\}$, where $\lambda \vdash n$ means $\lambda$ is an integer partition of $n$.

In the case that $G$ is a tree, the authors of \cite{gonzalez2024chromatic, GOT2} investigate the expansion $\textbf{X}_G$ in the star-basis, $\textbf{X}_G= \sum_{\lambda \vdash n} c_\lambda \mathfrak{st}_{\lambda} $. They focus on the \textit{leading partition} $\lead$ of $\mathbf{X}_G$, defined as the smallest partition in lexicographic order appearing in $\mathbf{X}_G$ with non-zero coefficient. The authors additionally provided formulas for the leading coefficient $c_\lead$, which encodes the degrees of the ``deep" vertices in $G$, and for the hook coefficients which determine the number of non-leaf (``internal") edges.

In this paper, we analyze the chromatic symmetric function of connected unicyclic graphs in the star-basis. We determine some structural features of these graphs from specific coefficients in their star-basis expansion. Our first main result, \Cref{thm:hook_coeff_general}, provides an explicit formula for the coefficients $c_\lambda$ in the star-basis expansion when $\lambda$ is a hook partition. As an immediate application, we show in \Cref{cor:cyc size from coeff} that the size of the cycle in any unicyclic graph can be directly determined from the coefficient $c_{(n)}$. 
Next, in \Cref{thm:leading}, we determine the leading partition $\lead$ for unicyclic graphs. Using $\lead$, we show in \Cref{cor.num.leaves} how to compute the number of leaves in the graph. We also give formulas for the leading coefficient $c_\lead$ for any unicyclic graph. 

With our results on the hook coefficients, leading partition, and the leading coefficient, we distinguish several classes of unicyclic graphs based on their CSF. In \Cref{cor:cyc size hooks}, we show that unicyclic graphs with different cycle sizes must have differing CSFs. 
Finally,  in \Cref{prop: cuttefish_distinguished} we prove that \textit{cuttlefish} graphs, which are a certain subclass of \textit{squids}, are distinguished from among all connected unicyclic graphs.

We also present computational data and examples of non-isomorphic graphs with a $4$-cycle that share the same CSF in \Cref{fig: n=12,fig:n=13}. This demonstrates that for bipartite graphs, the CSF alone does not necessarily distinguish non-isomorphic graphs. As a result, we establish that the tree isomorphism problem is independent of the bipartite property in trees.

Previous work has been done on computing the chromatic symmetric function of graphs with cycles. For example, see \cite{DVW17, martin2007distinguishingtreeschromaticsymmetric, ORELLANA20141, WW23}. We note that unicyclic graphs are not distinguished by their chromatic symmetric function \cite{ORELLANA20141, stanley_1995}. However, by analyzing the coefficients and partitions appearing in $\mathbf{X}_G$, we gain insight into when two unicyclic graphs can or cannot be isomorphic.

In \Cref{sec.general.graphs}, we give an explicit formula for the coefficient $c_{(n)}$ in $\mathbf{X}_{G}$ when $G$ is a connected bicyclic graph. This coefficient depends only on the sizes of the cycles in the graph and the number of edges shared by the two cycles. Additionally, we prove several results involving the sums of coefficients in $\mathbf{X}_G$ indexed by partitions of the same length, and use these results to prove that we can recover the sizes of the cycles in connected bicyclic graphs.

% %%%%%%%%%%%%%%%%%%%%%%%%%%%%%%%%

This paper is organized as follows. In \Cref{sec.def.prelims}, we provide the necessary background. In \Cref{sec.hook.coeff}, we describe the coefficients for hook partitions. In \Cref{sec.special.graphs}, we give formulas for the star-expansion of special families of graphs. In \Cref{sec.lead.partition,sec.lead.coeff}, we present results on the leading partition and leading coefficient. We use the result in \Cref{sec.lead.partition} to  prove that cuttlefish are distinguished (and reconstructible) from among all connected unicyclic graphs. In \Cref{sec.general.graphs}, we present results on bicyclic graphs. Finally in \Cref{sec: data}, we give examples and data.

\subsection*{Acknowledgments}This material is based on work supported by the National Science Foundation under Grant Number
DMS-1916439 while the authors were participating in the Mathematics Research Communities (MRC) 2024
Summer Conference at Beaver Hollow Conference Center in Java Center, New York. We extend our thanks
to the American Mathematical Society for their support. A. Bingham was partially supported by the grant DST/INT/RUS/RSF/P-41/2021 from the Department of Science \& Technology, Govt. of India, and FONDECYT-ANID grant 3250472. L. Johnston was partially supported by NSF grant DMS-2053350. R. Orellana was partially supported by NSF grant DMS-2153998. 
This material is also based upon work supported by the National Science Foundation, while A. Bingham, R. Orellana and J. Pan were in residence at the ICERM semester program ``Categorification and Computation in Algebraic Combinatorics'' in Fall 2025.
% The authors would also like to thank the anonymous referee for careful reading and useful comments, in particular for suggesting that we recover the sizes of the of the cycles in bicyclic graphs.

\section{Definitions and preliminaries}
\label{sec.def.prelims}
\subsection{Graph theory}
In this paper, we utilize standard definitions and fundamental concepts from graph theory.
For an overview of graph theory background, refer to, for example, \cite{west}.
For convenience and to establish the necessary notation for this paper, we recall some of the basic concepts here. 

%%%%%%%%%%%%%%%%%%%%%%%%%%%%
\subsubsection*{Graph basics}
A \defn{graph} $G$ is an ordered pair $(V,E)$, where $V=V(G)$ is the set of vertices of $G$ and $E=E(G)$ is its set of edges.
The number of vertices in $V$ is the \defn{order} of $G$, denoted by $|V|$.
The \defn{degree} of a vertex $v$ is the number of edges incident to it, and is denoted by $\deg(v)$.
The \defn{degree sequence} of a graph is the list of the degrees of all of its vertices, arranged in non-increasing order.

A graph $G$ is said to be \defn{simple} if $G$ contains no loops or multiple edges between two vertices.
All graphs in this paper are assumed to be simple, unless otherwise stated.
A \defn{path} in a graph $G$ is a sequence of distinct vertices $v_1, \ldots, v_n$ such that $v_iv_{i+1}$ is an edge in $E(G)$ for each $1 \le i \le n-1$.
A path on $n$ vertices in a graph is denoted by $P_n$.
A graph $G$ is \defn{connected} provided that there is a path between any two vertices in $V(G)$.
A \defn{cycle} is a graph with an equal number of vertices and edges whose vertices can be placed around a circle so that two vertices are adjacent if and only if they appear consecutively along the circle.

We denote by $C_n$ a cycle on $n$ vertices in a graph.
A graph is a \defn{tree} if there is a unique path between any two vertices.
Equivalently, a tree is a connected graph that contains no cycles.
A \defn{rooted tree} is a tree with a designated root vertex.
We call a tree comprised of a single vertex a \defn{trivial tree}, while any tree with more than one vertex is a \defn{non-trivial tree}.
A graph such that each of its connected components is a tree is called a \defn{forest}.

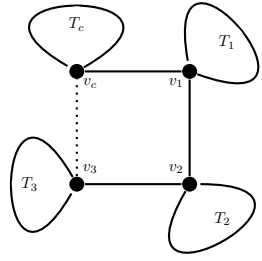
\begin{figure}[hbt!]
\centering
\begin{tikzpicture}[auto=center,every node/.style={circle, fill=black, scale=0.6}, style=thick, scale=0.5]
\node[label={[shift={(0.3,-0.8)}]$v_c$}] (Vc) at (0,0) {};
\node[label={[shift={(-0.3,-0.8)}]$v_1$}] (V1) at (3,0) {};
\node[label={[shift={(-0.3,-0.2)}]$v_2$}] (V2) at (3,-3) {};
\node[label={[shift={(0.3,-0.2)}]$v_3$}] (V3) at (0,-3) {};
\draw (Vc) -- (V1) -- (V2) -- (V3);
\draw[dotted] (V3)-- (Vc);
\draw[thick,-,shorten >=1pt] (V1) to [out=-20,in=100,loop,looseness=30,] node[below left,fill = white] {$T_1$}(V1);
\draw[thick,-,shorten >=1pt] (V2) to [out=-120,in=0,loop,looseness=30,] node[above left,above = 1.5pt, fill = white] {$T_2$}(V2);
\draw[thick,-,shorten >=1pt] (V3) to [out=120,in=240,loop,looseness=30,] node[right,fill = white] {$T_3$}(V3);
\draw[thick,-,shorten >=1pt] (Vc) to [out=30,in=150,loop,looseness=30,] node[below,fill = white] {$T_c$}(Vc);
\end{tikzpicture}
\caption{A connected unicyclic graph with a $c$-cycle and $c$ (potentially trivial) rooted trees.}
\label{fig.abstract.graph}
\end{figure}

A \defn{unicyclic graph} is a graph that contains exactly one cycle. Connected unicyclic graphs can also be characterized as simple graphs whose number of edges equals the number of vertices.
In this paper, we consider connected unicyclic graphs $G$ constructed from $c$ rooted trees $T_1, \ldots, T_c$ by adding edges between the roots of consecutive trees $T_i$ and $T_{i+1}$ for $1 \le i \le c-1$, 
and connecting the root of $T_c$ with the root of $T_1$, creating a cycle of size $c$ in the graph $G$ with vertices corresponding to the roots of the trees. We also call such a graph a \defn{$c$-unicyclic} graph. We use $r$ to denote the number of non-trivial rooted trees. See \Cref{fig.abstract.graph} for an illustration.

For example, the graph on the left of \Cref{eg.prelim} is a unicyclic graph constructed from $4$ rooted trees with roots $v_1,v_2,v_3$ and $v_4$ where $v_4$ is the root of a trivial tree. In this case, we have $r = 3$.
%%%%%%%%%%%%%%%%%%%%%%%%%%%%%%%
\subsubsection*{Vertex and edge types}
We now define specific types of vertices and edges in a graph that will be used throughout the paper.
An \defn{internal vertex} is a vertex with degree at least two.
A \defn{leaf} is a vertex of degree one. 
An internal vertex is a \defn{deep vertex} if it has no adjacent leaves. 
An \defn{internal edge} is an edge such that both of its endpoints are internal vertices, while a \defn{leaf-edge} is an edge such that at least one of its endpoints is a leaf. 
For a graph $G$, we denote by $I(G)$ the set of all internal edges of $G$. 
The connected components of the graph $G \setminus I(G)$ are referred to as the \defn{leaf components}.
For a connected unicyclic graph $G$ with $c$ rooted trees (as constructed above), we say vertex $v$ is a \defn{sprout} if it is a deep vertex on the $c$-cycle that is of degree at least three.

\begin{example}
Let $G$ be the graph in the left of \Cref{eg.prelim} on $14$ vertices. The internal vertices are $v_1,v_2,v_3,v_4,v_5,v_{11},v_{12}$, and the deep vertices are $v_3,v_4,v_{11}$. The only sprout is $v_3$.

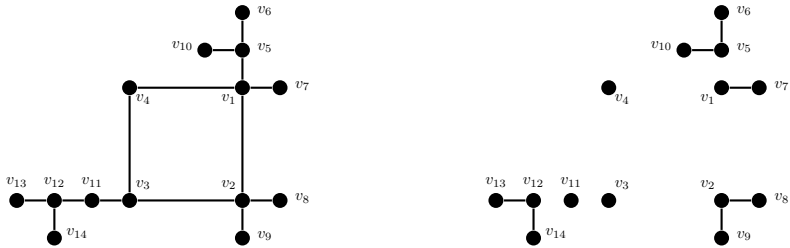
\begin{figure}[hbt!]
\centering
\begin{tikzpicture}[auto=center,every node/.style={circle, fill=black, scale=0.6}, style=thick, scale=0.5] 
\node[label={[shift={(0.3,-0.8)}]$v_4$}] (V4) at (0,0) {};
\node[label={[shift={(-0.3,-0.8)}]$v_1$}] (V1) at (3,0) {};
\node[label={[shift={(-0.3,-0.2)}]$v_2$}] (V2) at (3,-3) {};
\node[label={[shift={(0.3,-0.2)}]$v_3$}] (V3) at (0,-3) {};
\node[label={[shift={(0.5,-0.5)}]$v_5$}] (v5) at (3,1) {};
\node[label={[shift={(0.5,-0.5)}]$v_6$}] (v6) at (3,2) {};
\node[label={[shift={(-0.5,-0.5)}]$v_{10}$}] (v10) at (2,1) {};
\node[label={[shift={(0.5,-0.5)}]$v_8$}] (v8) at (4,-3) {};
\node[label={[shift={(0.5,-0.5)}]$v_9$}] (v9) at (3,-4) {};
\node[label={[shift={(0,-0.2)}]$v_{11}$}] (v11) at (-1,-3) {};
\node[label={[shift={(0,-0.2)}]$v_{12}$}] (v12) at (-2,-3) {};
\node[label={[shift={(0,-0.2)}]$v_{13}$}] (v13) at (-3,-3) {};
\node[label={[shift={(0.5,-0.5)}]$v_{14}$}] (v14) at (-2,-4) {};
\node[label={[shift={(0.5,-0.5)}]$v_7$}] (v7) at (4,0) {};
\draw (V1) -- (V2) -- (V3) -- (V4) -- (V1);
\draw(v7) -- (V1) -- (v5) -- (v6);
\draw(v10) -- (v5);
\draw(v9) -- (V2) -- (v8);
\draw(v12) -- (v11) -- (V3);
\draw(v14) -- (v12) -- (v13);
\end{tikzpicture}\qquad \qquad \qquad
\begin{tikzpicture}[auto=center,every node/.style={circle, fill=black, scale=0.6}, style=thick, scale=0.5] 
\node[label={[shift={(0.3,-0.8)}]$v_4$}] (V4) at (0,0) {};
\node[label={[shift={(-0.3,-0.8)}]$v_1$}] (V1) at (3,0) {};
\node[label={[shift={(-0.3,-0.2)}]$v_2$}] (V2) at (3,-3) {};
\node[label={[shift={(0.3,-0.2)}]$v_3$}] (V3) at (0,-3) {};
\node[label={[shift={(0.5,-0.5)}]$v_5$}] (v5) at (3,1) {};
\node[label={[shift={(0.5,-0.5)}]$v_6$}] (v6) at (3,2) {};
\node[label={[shift={(-0.5,-0.5)}]$v_{10}$}] (v10) at (2,1) {};
\node[label={[shift={(0.5,-0.5)}]$v_8$}] (v8) at (4,-3) {};
\node[label={[shift={(0.5,-0.5)}]$v_9$}] (v9) at (3,-4) {};
\node[label={[shift={(0,-0.2)}]$v_{11}$}] (v11) at (-1,-3) {};
\node[label={[shift={(0,-0.2)}]$v_{12}$}] (v12) at (-2,-3) {};
\node[label={[shift={(0,-0.2)}]$v_{13}$}] (v13) at (-3,-3) {};
\node[label={[shift={(0.5,-0.5)}]$v_{14}$}] (v14) at (-2,-4) {};
\node[label={[shift={(0.5,-0.5)}]$v_7$}] (v7) at (4,0) {};
\draw(v7) -- (V1);
\draw(v5) -- (v6);
\draw(v10) -- (v5);
\draw(v9) -- (V2) -- (v8);
\draw(v13) -- (v12) -- (v14);
\end{tikzpicture}
\caption{Left: A unicyclic graph $G$ on $14$ vertices; right: $G \setminus I(G)$.}
\label{eg.prelim}
\end{figure}
\end{example}

\subsection{Symmetric functions}\label{section:Sym}
A \defn{partition} is a sequence $\lambda=(\lambda_1, \lambda_2 \ldots , \lambda_k)$ of positive integers such that $\lambda_1 \geq \lambda_2 \geq \dots \geq \lambda_k$. We say that $\lambda$ is a partition of $n$, denoted $\lambda \vdash n$ or $|\lambda|=n$, if $\sum_i\lambda_i=n$. Each $\lambda_i$ is a \defn{part}, and the number of parts is the \defn{length} of $\lambda$, denoted $\ell(\lambda)$. If $\lambda$ has exactly $m_i$ parts equal to $i$ for $1 \leq i \leq n$, we may also denote $\lambda$ by $\lambda=(1^{m_1}2^{m_2}\ldots n^{m_n})$. 
A \defn{hook partition} of $n$ is a partition of the form $\la=(n-m_1, 1^{m_1})$, including the case where $m_1=0$.
The total ordering on partitions of $n$ that we use throughout the paper is the \defn{lexicographic order}:
$\mu \leq \lambda$ if $\mu=\lambda$ or $\mu_i=\lambda_i$ for $1 \leq i <j$ and $\mu_j < \lambda_j$ for some $1\leq j \leq \ell(\mu)$.
%%%%%%%%%%%%%%%%%%%%%%%%%
%%%%%%%%%%%%%%%%%%%%%%%%%%
\subsubsection*{Algebra of symmetric functions}
A formal power series $f$ in a countably infinite set of commuting variables $x_1, x_2, \dots$ is called a \defn{symmetric function} if it has bounded degree and is invariant under the action of any permutation $w$ of the positive integers on the indices of the variables, that is 
\[
f(x_1,x_2,\dots) = f(x_{w(1)},x_{w(2)},\dots)\,.
\]
Let $\Lambda^n$ be the vector space of all homogeneous symmetric functions of degree $n$.
The \defn{$r$-th elementary symmetric function} $e_r$ is defined by 
$$
e_r:=\sum_{i_1<i_2<\dots < i_r} x_{i_1}x_{i_2}\cdots x_{i_r}.
$$
For a partition $\lambda = (\lambda_1, \lambda_2, \dots, \lambda_k)$, the \defn{elementary symmetric function} corresponding to $\lambda$ is given by
$$
e_{\lambda} := e_{\lambda_1}e_{\lambda_2}\cdots e_{\lambda_k}
$$
forms a basis of $\Lambda^n$, called the \defn{elementary basis}.  Using the elementary symmetric functions, we can define the graded $\Q$-algebra of symmetric functions, denoted by $\Lambda$, as 
\[
\Lambda = \Lambda^0 \oplus \Lambda^1 \oplus \Lambda^2 \oplus \cdots\,,
 \text{ where } \Lambda^n := \text{Span}_{\mathbb{Q}}\{e_{\lambda} \ |\ \lambda \vdash n\} \text{ for } n \geq 1\,,
\]
and $\Lambda^0 :=\mathbb{Q}$.
%%%%%%%%%%%%%%%%%%%%%%%%%%
\subsubsection*{The chromatic symmetric function}
Let $G$ be a finite graph with vertex set $V=\{v_1, \ldots, v_n\}$. 
A \defn{proper coloring} of $G$ is a function $\kappa: V \to \mathbb{N}$ such that $\kappa(v_i) \neq \kappa(v_j)$ whenever $v_iv_j$ is an edge of $G$ with $1 \le i < j \le n$.
The \defn{chromatic symmetric function} (CSF) \cite{stanley_1995} of $G$ is defined by 
$$
{\bf X}_G = \sum_\kappa x_{\kappa(v_1)}x_{\kappa(v_2)}\dots x_{\kappa(v_n)},
$$
where the sum is over all proper colorings $\kappa$ of $G$.
The function ${\bf X}_G$ is a homogeneous symmetric function of degree $n$, where $n$ is the order of $G$. 
In the case that $G$ is a disjoint union of two subgraphs, $G=H_1 \sqcup H_2$, the chromatic symmetric function of $G$ is $\mathbf{X}_G=\mathbf{X}_{H_1}\cdot \mathbf{X}_{H_2}$.
In their work \cite{CVW16}, Cho and van Willigenburg give new bases of symmetric functions using the chromatic symmetric function and they proved the following result. 
\begin{theorem}[{\cite[Theorem 5]{CVW16}}]
\label{BasisTheorem}
    For any positive integer $k$, let $G_k$ denote a connected graph with $k$ vertices and let $\{G_k\}_{k \ge 1}$ be a family of such graphs. 
    Given a partition $\lambda \vdash n$, of length $\ell$, define $G_{\lambda}=G_{\lambda_1} \sqcup G_{\lambda_2} \sqcup \cdots \sqcup G_{\lambda_{\ell}}$.
    Then $\{{\bf X}_{G_{\lambda}} \ \mid \ \lambda \vdash n\}$   is a basis for $\Lambda^n$.
\end{theorem}

%%%%%%%%%%%%%%%%%%%%%%%%%%%
\subsubsection*{The star-basis}
In this paper, we focus on the star-basis of symmetric functions, which is constructed from the chromatic symmetric functions of star graphs. A \defn{star graph} on $k$ vertices, denoted by $\St_k$, is a tree with $k-1$ vertices of degree $1$ and one central vertex of degree $k-1$. We denote the chromatic symmetric function of $\St_k$ by $\mathfrak{st}_{k}:=\mathbf{X}_{\St_k}$.
For a partition $\lambda = (\lambda_1, \lambda_2, \dots, \lambda_{\ell})$, the chromatic symmetric function of the star forest $\St_{\lambda_1} \sqcup \St_{\lambda_2} \sqcup \cdots \sqcup \St_{\lambda_{\ell}}$ is denoted by 
$$
\mathfrak{st}_{\lambda}:=\mathfrak{st}_{\lambda_1}\mathfrak{st}_{\lambda_2} \cdots \mathfrak{st}_{\lambda_{\ell}}
\,.
$$
By \cref{BasisTheorem}, the set $\{\mathfrak{st}_{\lambda}\ \mid\ \lambda \vdash n\}$ is a basis, called the \defn{star-basis}, for the space $\Lambda^n$.
Hence, we may express ${\bf X}_G$ as:
$$ 
{\bf X}_G = \sum_{\lambda \vdash n} c_\lambda \mathfrak{st}_\lambda
\,.
$$

\begin{definition}\label{def:lead}
% Order the partitions of $n$ lexicographically in the order inherited from $\N^n$.
Given a graph $G$, we define the \defn{leading partition} to be the \emph{smallest} partition in lexicographic order such that $c_\lead\neq 0$, where $\mathbf{X}_G=\sum_\la c_\la \stf_\la$ is the star-expansion of $G$, which will be denoted $\lead(\mathbf{X}_G)$. When the context is clear, we will use $\lead$ for short.
\end{definition}
An algorithm is given in \cite{AMOZ23} for writing the chromatic symmetric function of a graph in the star-basis; this involves 
the deletion-near-contraction relations, which we recall next.

Given two sequences $a = (a_1,\dots,a_p)$ and $b = (b_1,\dots,b_q)$, let $(a,b) = a\cdot b = (a_1,\dots,a_p,b_1,\dots,b_q)$. Let \defn{$\sort()$} be the function that orders a sequence into a partition. The next lemma follows from the definitions of leading partition and lexicographic order. We will use it in the proofs in \Cref{sec.lead.partition,sec.lead.coeff}.

\begin{lemma}\label{lem.disjoint.lead}
Let $H_1$ and $H_2$ be two graphs. Then we have
$\lead(\mathbf{X}_{H_1\sqcup H_2}) = \sort(\lead(\mathbf{X}_{H_1})\cdot \lead(\mathbf{X}_{H_2}))$.
\end{lemma}

\subsection{Deletion-near-contraction relation}
Denote by $G\setminus e$ the graph obtained from $G$ by deleting edge $e$. The \defn{leaf-contraction} operation is denoted $G\odot e$ and yields the graph obtained by contracting edge $e$ and attaching a new leaf to the vertex that results from this contraction, connected by a new leaf-edge $\ell_e$ . The \defn{dot-contraction} $(G\odot e) \setminus \ell_e$ amounts to contracting edge $e$ and then adding an isolated vertex. 

We note that if $G$ contains a 3-cycle and $e$ is an edge on the cycle, then $G\odot e$ and $(G\odot e) \setminus \ell_e$ will produce graphs with multiple edges. To ensure that all graphs are simple, we will delete any duplicate edges in these cases. 
\begin{proposition}
[\cite{AMOZ23}][The deletion-near-contraction (DNC) relation]
\label{DNC.relation}
Let $G$ be a simple graph and $e$ be any edge in $G$. Then 
$$
\mathbf{X}_G= \mathbf{X}_{G\setminus e} - \mathbf{X}_{(G \odot e)\setminus {\ell_e}} + \mathbf{X}_{G \odot e}. $$
\end{proposition}

Observe that if $e$ is a leaf-edge, then the DNC relation does not provide any new information for $\mathbf{X}_G$, as $G\setminus e\cong(G \odot e)\setminus {\ell_e}$ and $G \odot e\cong G$. Thus, we will only apply the DNC relation to internal edges of $G$. Furthermore if $e$ is an internal edge, then the graphs $G\setminus e,\  (G \odot e)\setminus {\ell_e},$ and $G \odot e$ each have fewer internal edges than $G$. Note that a simple graph with no internal edges is a star forest. Therefore, by repeatedly applying the DNC relation to internal edges, we can express $\mathbf{X}_G$ as a linear combination of $\mathfrak{st}_{\lambda}$, where $\lambda$ is a partition of the number of vertices of $G$.

We now recall the star-expansion algorithm introduced in \cite{AMOZ23}, which formalizes the recursive construction of $\mathbf{X}_G$ in the star-basis described above.\vspace{1em}

\begin{algorithm}[H]
\textbf{Input:} A simple graph $G$.\\
\textbf{Initialization:} Let $\mathcal{T}$ be a tree with root labeled $G$ and no edges.\\
\textbf{Iteration:} If $H$ is a leaf $\mathcal{T}$ and $H$ has some internal edge $e$, then extend $\mathcal{T}$ by adding three children to $H$ labeled $H\setminus e, (H \odot e)\setminus {\ell_e}$, and $H \odot e$. Label the three new edges with $+$ or $-$ according to the corresponding coefficient in the DNC relation. The algorithm terminates when all leaves in $\mathcal{T}$ contain no internal edges.\\
\textbf{Output:} A rooted tree $\mathcal{T}(G)$ with leaves labeled by star forests.
\caption{Star-expansion}
\end{algorithm}

We call the output of the star-expansion algorithm the \defn{DNC-tree}. Using the DNC-tree, the authors in \cite{AMOZ23} demonstrated how to write the star-expansion of $\mathbf{X}_G$. We first introduce some notation. Let $\iota(H)$ and $\iota(G)$ denote the number of isolated vertices of the graphs $H$ and $G$, respectively. Let $\lambda(H)$ be the partition whose parts correspond to the orders of the connected components of $H$. 

\begin{theorem}[\cite{AMOZ23}]
For every simple graph $G$, let $\mathcal{T}(G)$ be a DNC-tree obtained from the star-expansion algorithm, and let $L(\mathcal{T}(G))$ be the multiset of leaf labels of $\mathcal{T}(G)$. Then
\[
\mathbf{X}_G = \sum_{H \in L(\mathcal{T}(G))} (-1)^{\iota(H)-
\iota(G)}\mathfrak{st}_{\lambda(H)}.
\]
In addition, no cancellations occur in the computation; that is, for any partition $\lambda$, all terms $\mathfrak{st}_\lambda$ appear with the same sign. 
\end{theorem}

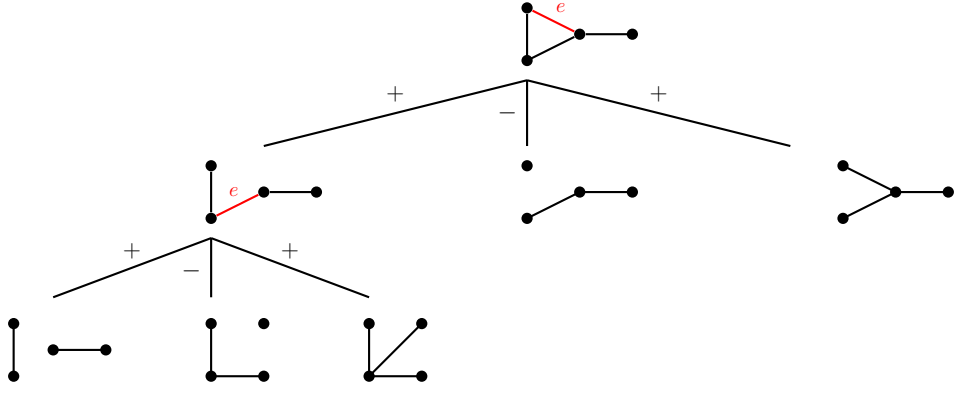
\begin{figure}[hbt!]\label{Ex: DNC-tree}
    \centering
    \begin{tikzpicture}[auto=center,every node/.style={circle, fill=black, scale=0.45}, style=thick, scale=0.35] \label{tikz:dnc-tree}
    %%%%% TREE 1 %%%%%
     \node (A1) at (0,1) {};
    \filldraw[black] (0, -1) coordinate (A2) circle (4pt) node{};
     \node (A3) at (2,0) {};
    \filldraw[black] (4, 0) coordinate (A4) circle (4pt) node{};

    \draw(A1) -- (A2);
    \draw(A2) -- (A3);
    \draw[color=red](A3) -- (A1) node[above,right = 8pt,fill=white] {\huge{$e$}};
    \draw(A3) -- (A4);

     %%%%% TREE 2 %%%%%
     \node (A1) at (-12,-5) {};
     \node (A2) at (-12, -7) {};
     \node (A3) at (-10,-6) {};
    \filldraw[black] (-8, -6) coordinate (A4) circle (4pt) node{};

    \draw(A1) -- (A2);
    \draw[color=red](A2) -- (A3) node[above,right = -16pt,fill=white] {\huge{$e$}};
    \draw(A3) -- (A4);
    
    %%%%% TREE 3 %%%%%
    \filldraw[black] (0,-5) coordinate (A1) circle (4pt) node{};
    \filldraw[black] (0, -7) coordinate (A2) circle (4pt) node{};
    \filldraw[black] (2, -6) coordinate (A3) circle (4pt) node{};
    \filldraw[black] (4, -6) coordinate (A4) circle (4pt) node{};
    
    \draw (A2) -- (A3);
    \draw(A3) -- (A4);

     %%%%% TREE 4 %%%%%
    \filldraw[black] (12,-5) coordinate (A1) circle (4pt) node{};
    \filldraw[black] (12, -7) coordinate (A2) circle (4pt) node{};
    \filldraw[black] (14, -6) coordinate (A3) circle (4pt) node{};
    \filldraw[black] (16, -6) coordinate (A4) circle (4pt) node{};

    \draw(A1) -- (A3);
    \draw (A2) -- (A3);
    \draw(A3) -- (A4);

    %%%%% TREE 5 %%%%%
    \filldraw[black] (-19.5,-11) coordinate (A1) circle (4pt) node{};
    \filldraw[black] (-19.5, -13) coordinate (A2) circle (4pt) node{};
    \filldraw[black](-18,-12) coordinate (A3) circle (4pt) node{};
    \filldraw[black] (-16, -12) coordinate (A4) circle (4pt) node{};

    \draw(A1) -- (A2);
    \draw(A3) -- (A4);
    
    %%%%% TREE 6 %%%%%
    \filldraw[black] (-12,-11) coordinate (A1) circle (4pt) node{};
    \filldraw[black] (-12, -13) coordinate (A2) circle (4pt) node{};
    \filldraw[black](-10,-13) coordinate (A3) circle (4pt) node{};
    \filldraw[black] (-10, -11) coordinate (A4) circle (4pt) node{};

    \draw(A1) -- (A2);
    \draw(A2) -- (A3);

    %%%%% TREE 7 %%%%%
    \filldraw[black] (-6,-11) coordinate (A1) circle (4pt) node{};
    \filldraw[black] (-6, -13) coordinate (A2) circle (4pt) node{};
    \filldraw[black](-4,-13) coordinate (A3) circle (4pt) node{};
    \filldraw[black] (-4, -11) coordinate (A4) circle (4pt) node{};

    \draw(A1) -- (A2);
    \draw(A2) -- (A3);
    \draw(A2) -- (A4);

    %%%%% DNC-tree EDGES %%%%
    \draw(0,-1.75) -- (-10,-4.25);
    \draw(0,-1.75) -- (0,-4.25);
    \draw(0,-1.75) -- (10,-4.25);

    \node[fill=none] at (-5, -2.25) {\Huge $+$};
    \node[fill=none] at (-.75, -3) {\Huge $-$};
    \node[fill=none] at (5, -2.25) {\Huge $+$};

    \draw(-12,-7.75) -- (-18,-10);
    \draw(-12,-7.75) -- (-12,-10);
    \draw(-12,-7.75) -- (-6,-10);

    \node[fill=none] at (-15, -8.25) {\Huge $+$};
    \node[fill=none] at (-12.75, -9) {\Huge $-$};
    \node[fill=none] at (-9, -8.25) {\Huge $+$};
    \end{tikzpicture}
    
    \caption{The DNC-tree $\mathcal{T}(G)$.}
    \label{fig:dnc-tree}
\end{figure}

\begin{example}
    We give an example of the star-expansion algorithm in Figure~\ref{fig:dnc-tree} for the graph $G$ at the top of the DNC-tree. The red edges indicate to which edge the DNC relation is applied. The resulting CSF in the star-basis is \[\mathbf{X}_G=2\mf{st}_{(4)}-2\mf{st}_{(3,1)}+\mf{st}_{(2,2)}.\]
\end{example}

\section{Coefficients for hook partitions}

\label{sec.hook.coeff}
We give the coefficients for hook partitions for all connected unicyclic graphs. First, we give a known fact about the coefficients of hook partitions for trees. 
\begin{proposition}[\cite{gonzalez2024chromatic}, Proposition 3.8]
\label{prop:treehook}
    Let $T$ be an $n$-vertex tree and $k$ be the number of internal edges of $T$. If $\mathbf{X}_T= \sum_{\lambda \vdash n}c_\lambda \mathfrak{st}_\lambda$, then 
    \[ c_{(n-m_1,1^{m_1})}=(-1)^{m_1} \binom{k}{m_1}.\] 
\end{proposition}

Note that a unicyclic graph with cycle size $c$ always has at least $c$ internal edges. The following lemma is the base case for Theorem~\ref{thm:hook_coeff_general}, which is proved by induction on the size of the cycle. 

\begin{lemma}
\label{lemma:3_cycle_hooks}
    Let $G$ be a connected unicyclic graph on $n$ vertices with one 3-cycle, let $\lambda= (n-m_1,1^{m_1})$ be a hook partition  (note $m_1\neq n-1$), and let $k$ be the number of internal edges of $G$. Let $0 \leq r \leq 3$ be the number of (non-trivial) rooted trees on the $3$-cycle. If $\mathbf{X}_G= \sum_{\lambda \vdash n}c_\lambda \mathfrak{st}_\lambda$ then 
        \[c_{(n-m_1,1^{m_1})}=  (-1)^{m_1} \left[(r-1)\binom{k-2}{m_1-1}+2 \binom{k-2}{m_1}\right].
        \]
\end{lemma}
\begin{proof}
We prove each case for $r$ separately. 
    \begin{enumerate}
    \setcounter{enumi}{-1}
        \item Observe that when $r=0$, the graph $G$ is the cycle graph on 3 vertices, $C_3$, and $\mathbf{X}_{C_3}=2 \mathfrak{st}_{(3)}-1 \mathfrak{st}_{(2,1)}$. This agrees with the formula since $k=3$. 
        \item Suppose $r=1$, and let $e$ be an edge on the cycle that is adjacent to the single rooted tree. The graph $G\setminus e$ is a tree with $k-2$ internal edges, so the coefficient $c_{(n-m_1,1^{m_1})}$ in $\mathbf{X}_{G\setminus e}$ is $(-1)^{m_1}\binom{k-2}{m_1}$ by Proposition~\ref{prop:treehook}. The graph $(G \odot e)\setminus {\ell_e}$ is the union of an isolated vertex and a tree with $k-3$ internal edges, so the coefficient $c_{(n-m_1,1^{m_1})}$ in $\mathbf{X}_{(G \odot e)\setminus {\ell_e}}$ is $(-1)^{m_1-1}\binom{k-3}{m_1-1}$. The graph $G \odot e$ is a tree with $k-3$ internal edges, so the $c_{(n-m_1,1^{m_1})}$ in $\mathbf{X}_{(G \odot e)}$ is $(-1)^{m_1}\binom{k-3}{m_1}$. 
        \cref{DNC.relation} and Pascal's rule imply that in $\mathbf{X}_G$, 
        $c_{(n-m_1,1^{m_1})}=(-1)^{m_1}\left[ \binom{k-2}{m_1}+\binom{k-3}{m_1-1}+\binom{k-3}{m_1}\right]= (-1)^{m_1} \left[\binom{k-2}{m_1}+\binom{k-2}{m_1}\right]=(-1)^{m_1} \left[2\binom{k-2}{m_1}\right]$. 
        \item Suppose $r=2$, and let $e$ be an edge on the cycle that is adjacent to only a single rooted tree. Following a similar argument as above, the graph $G\setminus e$ is a tree with $k-2$ internal edges, so the coefficient $c_{(n-m_1,1^{m_1})}$ in $\mathbf{X}_{G\setminus e}$ is $(-1)^{m_1}\binom{k-2}{m_1}$. The graph $(G \odot e)\setminus {\ell_e}$ is the union of an isolated vertex and a tree with $k-2$ internal edges, so $c_{(n-m_1,1^{m_1})}$ in $\mathbf{X}_{(G \odot e)\setminus {\ell_e}}$ is $(-1)^{m_1-1}\binom{k-2}{m_1-1}$. The graph $G \odot e$ is a tree with $k-2$ internal edges, so $c_{(n-m_1,1^{m_1})}$ in $\mathbf{X}_{(G \odot e)}$ is $(-1)^{m_1}\binom{k-2}{m_1}$. Thus in $\mathbf{X}_G$, $c_{(n-m_1,1^{m_1})}=(-1)^{m_1}\left[ \binom{k-2}{m_1}+\binom{k-2}{m_1-1}+\binom{k-2}{m_1}\right]=(-1)^{m_1}\left[ \binom{k-2}{m_1-1}+2\binom{k-2}{m_1}\right]$.
        \item Let $r=3$, and let $e$ be any edge on the cycle. The graph $G\setminus e$ is a tree with $k-1$ internal edges, so the coefficient $c_{(n-m_1,1^{m_1})}$ in $\mathbf{X}_{G\setminus e}$ is $(-1)^{m_1}\binom{k-1}{m_1}$. The graph $(G \odot e)\setminus {\ell_e}$ is the union of an isolated vertex and a tree with $k-2$ internal edges, so $c_{(n-m_1,1^{m_1})}$ in $\mathbf{X}_{(G \odot e)\setminus {\ell_e}}$ is $(-1)^{m_1-1}\binom{k-2}{m_1-1}$. The graph $G \odot e$ is a tree with $k-2$ internal edges, so $c_{(n-m_1,1^{m_1})}$ in $\mathbf{X}_{(G \odot e)}$ is $(-1)^{m_1}\binom{k-2}{m_1}$. Thus in $\mathbf{X}_G$, 
        \begin{align*}
            c_{(n-m_1,1^{m_1})} =& (-1)^{m_1}\left[ \binom{k-1}{m_1}+\binom{k-2}{m_1-1}+\binom{k-2}{m_1}\right]\\
            =& (-1)^{m_1} \left[ 2\binom{k-2}{m_1-1}+2\binom{k-2}{m_1}\right].
        \end{align*}   
    \end{enumerate}
\end{proof}

The following theorem, which we will prove by induction on the cycle size, covers the general case of hook coefficients for connected unicyclic graphs, where \Cref{lemma:3_cycle_hooks} is the base case.

\begin{theorem}
\label{thm:hook_coeff_general}
Let $G$ be a connected unicyclic graph on $n$ vertices with a single cycle of size $c\geq 3$, let $\lambda= (n-m_1,1^{m_1})$ be a hook partition, $k$ be the number of internal edges of $G$, and $r\geq 0$ be the number of non-trivial rooted trees attached to the cycle in $G$. If $\mathbf{X}_G= \sum_{\lambda \vdash n}c_\lambda \mathfrak{st}_\lambda$, then 
\begin{equation}\label{eq:hook_coeff}
c_{(n-m_1,1^{m_1})}= (-1)^{m_1} \left[(r-1)\binom{k-2}{m_1-1} + (c-1) \binom{k-2}{m_1} \right].
\end{equation}
\end{theorem}

\begin{proof} We will split into cases based on the number of rooted trees.

\textbf{Case 1 ($r\geq 1$):} We proceed by inducting on the size of the cycle $c$. The base case when $c=3$ is Lemma~\ref{lemma:3_cycle_hooks}. Now suppose that the claim is true for any unicyclic graph with a cycle of size $c=z\geq 3$. Let $G$ be a connected graph on $n$ vertices with a single cycle of size $z+1$, $r\geq 1$ rooted trees, and $k$ internal edges. We show that the coefficient $c_{(n-m_1,1^{m_1})}$ in $\mathbf{X}_G$ is $(-1)^{m_1}\left[ (r-1)\binom{k-2}{m_1-1}+z\binom{k-2}{m_1}\right]$. Since $r\geq 1$ and the rooted trees can be placed anywhere on the cycle, there is either an edge on the cycle that is adjacent to only a single rooted tree, or there exists no such edge and hence every edge of the cycle is adjacent to two rooted trees. We examine each case separately.\vspace{1em}

\textbf{Case 1a:} Suppose $e$ is an edge on the cycle that is adjacent to only a single rooted tree. $G \setminus e$ is a tree with $k-2$ internal edges, so the coefficient $c_{(n-m_1,1^{m_1})}$ in $\mathbf{X}_{G\setminus e}$ is $(-1)^{m_1} \binom{k-2}{m_1}$. Next, $(G \odot e)\setminus {\ell_e}$ is the union of an isolated vertex with a unicyclic graph with cycle size $z$, $r$ rooted trees, and $k-1$ internal edges. By our inductive hypothesis, the coefficient $c_{(n-m_1,1^{m_1})}$ in $\mathbf{X}_{(G \odot e)\setminus {\ell_e}}$ is $(-1)^{m_1-1}\left[ (r-1)\binom{k-3}{m_1-2}+(z-1)\binom{k-3}{m_1-1}\right]$. Finally, $G \odot e$ is a unicyclic graph with cycle size $z$, $r$ rooted trees, and $k-1$ internal edges. Again by our inductive hypothesis, $c_{(n-m_1,1^{m_1})}$ in $\mathbf{X}_{G \odot e}$ is $(-1)^{m_1}\left[(r-1)\binom{k-3}{m_1-1} + (z-1)\binom{k-3}{m_1}\right]$. Thus by the DNC relation and Pascal's rule, the coefficient of $c_{(n-m_1,1^{m_1})}$ in $\mathbf{X}_G$ is
\begin{align*}
    c_{(n-m_1,1^{m_1})} =& (-1)^{m_1} \left[ \binom{k-2}{m_1}+ (r-1)\left[\binom{k-3}{m_1-2}+\binom{k-3}{m_1-1}\right]+(z-1)\left[\binom{k-3}{m_1-1}+\binom{k-3}{m_1}\right]\right]\\
    =& (-1)^{m_1} \left[ \binom{k-2}{m_1}+ (r-1)\binom{k-2}{m_1-1}+(z-1)\binom{k-2}{m_1}\right]\\
    =& (-1)^{m_1}\left[ (r-1)\binom{k-2}{m_1-1}+z\binom{k-2}{m_1}\right].
\end{align*}

\textbf{Case 1b:} Suppose $e$ is an edge on the cycle that is adjacent to two rooted trees. So $r\geq 2$. Following a similar argument as above, $G \setminus e$ is a tree with $k-1$ internal edges, so the coefficient $c_{(n-m_1,1^{m_1})}$ in $\mathbf{X}_{G\setminus e}$ is $(-1)^{m_1} \binom{k-1}{m_1}$. Next, $(G \odot e)\setminus {\ell_e}$ is the union of an isolated vertex with a unicyclic graph with cycle size $z$, $r-1$ rooted trees, and $k-1$ internal edges. By our inductive hypothesis, the coefficient $c_{(n-m_1,1^{m_1})}$ in $\mathbf{X}_{(G \odot e)\setminus {\ell_e}}$ is $(-1)^{m_1-1}\left[ (r-2)\binom{k-3}{m_1-2}+(z-1)\binom{k-3}{m_1-1}\right]$. Finally, $G \odot e$ is a unicyclic graph with cycle size $z$, $r-1$ rooted trees, and $k-1$ internal edges. By our inductive hypothesis, $c_{(n-m_1,1^{m_1})}$ in $\mathbf{X}_{G \odot e}$ is $(-1)^{m_1}\left[(r-2)\binom{k-3}{m_1-1} + (z-1)\binom{k-3}{m_1}\right]$. Thus the coefficient $c_{(n-m_1,1^{m_1})}$ in $\mathbf{X}_G$ is
\begin{align*}
    c_{(n-m_1,1^{m_1})} =& (-1)^{m_1} \left[ \binom{k-1}{m_1}+ (r-2)\left[\binom{k-3}{m_1-2}+\binom{k-3}{m_1-1}\right]+(z-1)\left[\binom{k-3}{m_1-1}+\binom{k-3}{m_1}\right]\right]\\
    =& (-1)^{m_1} \left[ \binom{k-1}{m_1}+ (r-2)\binom{k-2}{m_1-1}+(z-1)\binom{k-2}{m_1}\right]\\
    =& (-1)^{m_1} \left[ \binom{k-2}{m_1}+ \binom{k-2}{m_1-1}+ (r-2)\binom{k-2}{m_1-1}+(z-1)\binom{k-2}{m_1}\right]\\
    =& (-1)^{m_1}\left[ (r-1)\binom{k-2}{m_1-1}+z\binom{k-2}{m_1}\right].
\end{align*}

\textbf{Case 2 ($r= 0$)}: We again induct on the size of the cycle $c$. Suppose that for any cycle graph $C_z$ on $z$ vertices where $z\geq 3$, the hook coefficients are \[c_{(z-m_1,1^{m_1})}=(-1)^{m_1}\left[-\binom{z-2}{m_1-1}+(z-1)\binom{z-2}{m_1}\right].\] This follows from $r=0$ and $n=k=z$ in this case. Let $G=C_{z+1}$. We show that the hook coefficients in $C_{z+1}$ are \[c_{(z+1-m_1,1^{m_1})}=(-1)^{m_1}\left[-\binom{z-1}{m_1-1}+z\binom{z-1}{m_1}\right].\]
Let $e$ be any edge on the cycle. Then $G \setminus e$ is a tree with $z-2$ internal edges, so $c_{(z+1-m_1,1^{m_1})}$ in $\mathbf{X}_{G\setminus e}$ is $(-1)^{m_1} \binom{z-2}{m_1}$.
Now, $(G \odot e)\setminus {\ell_e}$ is the union of an isolated vertex with the cycle graph on $z$ vertices. By our inductive hypothesis, the coefficient $c_{(z+1-m_1,1^{m_1})}$ in $\mathbf{X}_{(G \odot e)\setminus {\ell_e}}$ is $(-1)^{m_1-1}\left[ -\binom{z-2}{m_1-2}+(z-1)\binom{z-2}{m_1-1}\right]$. Finally, $G \odot e$ is a unicyclic graph with cycle size $z$, $1$ rooted tree, and $z$ internal edges. By case 1, $c_{(z+1-m_1,1^{m_1})}$ in $\mathbf{X}_{G \odot e}$ is $(-1)^{m_1}\left[(z-1)\binom{z-2}{m_1}\right]$. Thus the coefficient $c_{(z+1-m_1,1^{m_1})}$ in $\mathbf{X}_G$ is
\begin{align*}
    c_{(n-m_1,1^{m_1})} =& (-1)^{m_1} \left[ \binom{z-2}{m_1}-\binom{z-2}{m_1-2}+(z-1)\binom{z-2}{m_1-1}+(z-1)\binom{z-2}{m_1}\right]\\
    =& (-1)^{m_1} \left[ \binom{z-2}{m_1}+\binom{z-2}{m_1-1}-\binom{z-1}{m_1-1}+(z-1)\binom{z-1}{m_1}\right]\\
    =& (-1)^{m_1}\left[-\binom{z-1}{m_1-1}+z\binom{z-1}{m_1}\right].
\end{align*}
\end{proof}

Consequently, given any connected unicyclic graph $G$, we can determine the coefficients for the elements in the star-basis expansion corresponding to hook partitions. We immediately obtain the following corollaries.

\begin{corollary}\label{cor:cyc size from coeff}
    Let $G$ be a connected unicyclic graph of order $n$ with a single cycle of size $c \geq 3$. If $\mathbf{X}_G = \sum_{\lambda \vdash n} c_\lambda \mathfrak{st}_\lambda$, then
    \[ c_{(n)}=c-1.\]
\end{corollary}
\begin{corollary} \label{cor:cyc size hooks}
    If $G_1$ and $G_2$ are connected unicyclic graphs with cycles of different sizes, then $\mathbf{X}_{G_1} \neq \mathbf{X}_{G_2}$.
\end{corollary}
In particular, the previous corollary shows that the star-expansion of the chromatic symmetric function of a unicyclic graph can be used to distinguish the cycle $C_n$ from other connected unicyclic graphs of order $n$. The cycle will be the only connected unicyclic graph with $c_{(n)}=n-1$; alternatively, it will be the only connected unicyclic graph for which $c_{(2,1^{n-2})}\neq 0$, since no other connected unicyclic graph with $n$ vertices satisfies $k=n$.

Next, we show how the longest hook partition encodes the relationship between the number of internal edges and rooted trees.
\begin{proposition}
 \label{prop:k-and-r-hooks}
    Let G be a connected unicyclic graph of order $n$, $k$ internal edges, $r$ rooted trees, and a cycle of size $c \geq 3$. If $\mathbf{X}_G = \sum_{\lambda \vdash n} c_\lambda \mathfrak{st}_\lambda$ and $\lambda^\star=(n-m_1,1^{m_1})$ is the longest hook partition  such that $c_{\lambda^\star} \neq 0$, then 
    \[ k=\begin{cases}
        m_1+2 & \text{if } r=0 \text{ or } 1, \\
        m_1+1 & \text{if } r\geq 2 .
    \end{cases}\]
    Furthermore, the coefficient $c_{\lambda^\star}$ corresponding to the longest hook partition is 
     \[ c_{\lambda^\star}=
    \begin{cases}
        (-1)^{k-2} & \text{if } r=0, \\
        (-1)^{k-2}(c-1) & \text{if } r=1, \\
        (-1)^{k-1}(r-1) & \text{if } r\geq 2 .
    \end{cases}
    \]
\end{proposition}

\begin{proof}
    We consider the cases when $r=0$, $r=1$, and $r \geq 2$ separately. 
    First consider when $r=0$, that is, $G$ is cycle graph $C_n$, and $n=k=c$. Then Equation~(\ref{eq:hook_coeff}) reduces to
    \[ c_{(n-m_1,1^{m_1})}= (-1)^{m_1}\left[-\binom{k-2}{m_1-1}+(k-1)\binom{k-2}{m_1}\right].\]
    The largest $m_1$ such that $c_{(n-m_1,1^{m_1})} \neq 0$ is $m_1=k-2$, since the next hook longer than the partition $(n-k+2,1^{k-2})=(2,1^{k-2})$ is $(1^k)$, which is when $m_1=k$ and hence $c_{(1^k)}=0$.
    Next, consider when $r=1$. Then Equation~(\ref{eq:hook_coeff}) reduces to \[ c_{(n-m_1,1^{m_1})}= (-1)^{m_1}\left[(c-1)\binom{k-2}{m_1}\right]\] and we again see that the largest $m_1$ such that $c_{(n-m_1,1^{m_1})} \neq 0$ is $m_1=k-2$. Finally, if $r\geq 2$, the largest $m_1$ such that $c_{(n-m_1,1^{m_1})} \neq 0$ is $m_1=k-1$, as can be seen by examining Equation~(\ref{eq:hook_coeff}). 
\end{proof}

\begin{example} Let $G$ be the graph in Figure \ref{fig:hook-coeff}. We have that $\mathbf{X}_G=\textcolor{red}{3}\, \mf{st}_{(8)}\textcolor{red}{-10}\,\mf{st}_{(7,1)} +4\,\mf{st}_{(6,2)}+\textcolor{red}{12}\,\mf{st}_{(6,1,1)}+2\mf{st}_{(5,3)}-10\,\mf{st}_{(5,2,1)}\textcolor{red}{-6}\,\mf{st}_{(5,1,1,1)}+1\,\mf{st}_{(4,4)}-5\,\mf{st}_{(4,3,1)}+2\,\mf{st}_{(4,2,2)}+8\,\mf{st}_{(4,2,1,1)}+\textcolor{red}{1}\,\mf{st}_{(4,1,1,1,1)}+2\,\mf{st}_{(3,3,2)}+2\,\mf{st}_{(3,3,1,1)}-4\,\mf{st}_{(3,2,2,1)}-2\,\mf{st}_{(3,2,1,1,1)}+1\,\mf{st}_{(2,2,2,1,1)}$. The hook coefficients are highlighted in red and match Equation \ref{eq:hook_coeff} as $n=8, c=4, k=5$, and $r=2$.\par
Additionally, the longest hook partition is $\lambda^\star=(4,1,1,1,1)$ and the corresponding coefficient is $c_{\lambda^\star}=(-1)^{k-1}(r-1)=1$.

 \begin{figure}[hbt!]
    \centering
    \begin{tikzpicture}[auto=center,every node/.style={circle, fill=black, scale=0.45}, style=thick, scale=0.6] 
    %%%%% Draw vertices %%%%%
    \filldraw[black] (0, 1) coordinate (A0) circle (4pt) node{};
    \filldraw[black] (0, 0) coordinate (A1) circle (4pt) node{};
    \filldraw[black] (1, 0) coordinate (A2) circle (4pt) node{};
    \filldraw[black] (1, 1) coordinate (A3) circle (4pt) node{};
    \filldraw[black] (-1, 0) coordinate (A4) circle (4pt) node{};
    \filldraw[black] (1, 2) coordinate (A5) circle (4pt) node{};
    \filldraw[black] (2, 1) coordinate (A6) circle (4pt) node{};
    \filldraw[black] (3, 1) coordinate (A7) circle (4pt) node{};    
    
    %%%%% Draw edges %%%%%
    \draw(A1) -- (A2) -- (A3) -- (A0) -- (A1);
    \draw (A1) -- (A4);
    \draw (A3) -- (A5);
    \draw (A3) -- (A6) -- (A7);
    \end{tikzpicture}
     \caption{A unicyclic graph $G$ on $8$ vertices.}
     \label{fig:hook-coeff}
\end{figure}    
\end{example}

\begin{corollary} \label{cor:hook k r dist}
    Let $G_1$ and $G_2$ be connected unicyclic graphs on $n$ vertices with $k_1$ and $k_2$ internal edges and $r_1$ and $r_2$ rooted tree respectively. If $k_1 = k_2$ and $r_1 \neq r_2$ then $\mathbf{X}_{G_1} \neq \mathbf{X}_{G_2}$.
\end{corollary}

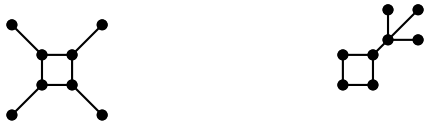
\begin{figure}[hbt!]
    \centering
    \begin{tikzpicture}[auto=center,every node/.style={circle, fill=black, scale=0.45}, style=thick, scale=0.5] 
        %%%%% TREE 1 %%%%%
    %%%%% Draw vertices %%%%%
    \filldraw[black] (0, 1) coordinate (A0) circle (4pt) node{};
    \filldraw[black] (0, 0) coordinate (A1) circle (4pt) node{};
    \filldraw[black] (1, 0) coordinate (A2) circle (4pt) node{};
    \filldraw[black] (1, 1) coordinate (A3) circle (4pt) node{};
    \filldraw[black] (2, -1) coordinate (A4) circle (4pt) node{};
    \filldraw[black] (2, 2) coordinate (A5) circle (4pt) node{};
    \filldraw[black] (-1,-1) coordinate (A6) circle (4pt) node{};
    \filldraw[black] (-1, 2) coordinate (A7) circle (4pt) node{};   
    
    %%%%% Draw edges %%%%%
    \draw(A1) -- (A2) -- (A3) -- (A0) -- (A1);
    \draw(A3) -- (A5);
    \draw(A1) -- (A6);
    \draw(A0) -- (A7);
    \draw(A2) -- (A4);

        %%%%% TREE 2 %%%%%
    %%%%% Draw vertices %%%%%
    \filldraw[black, ] (10, 1) coordinate (A0) circle (4pt) node{};
    \filldraw[black] (10, 0) coordinate (A1) circle (4pt) node{};
    \filldraw[black] (11, 0) coordinate (A2) circle (4pt) node{};
    \filldraw[black] (11, 1) coordinate (A3) circle (4pt) node{};
    \filldraw[black] (11.5, 1.5) coordinate (A4) circle (4pt) node{};
    \filldraw[black] (11.5, 2.5) coordinate (A5) circle (4pt) node{};
    \filldraw[black] (12.5, 1.5) coordinate (A6) circle (4pt) node{};    
    \filldraw[black] (12.5, 2.5) coordinate (A7) circle (4pt) node{};

    %%%%% Draw edges %%%%%
    \draw(A1) -- (A2)--(A3)--(A0)--(A1);
    \draw(A3)--(A4)--(A5);
    \draw(A7)--(A4)--(A6);

    \end{tikzpicture}
     \caption{Two graphs with the same hook coefficients: $c_{(8)}=3$, $c_{(7,1)}=-9$, $c_{(6,1,1)}=9$, $c_{(5,1^3)}=-3$, and $c_\la=0$ for all other hooks $\la$.}
     \label{fig:hooks-kr}
\end{figure}

    The coefficient $c_{(n-1,1)}$ can be used to write a quadratic equation for $r$ or $k$, yielding two possible solutions. We note however, that hook coefficients alone are insufficient for recovering these values exactly. As an example see Figure~\ref{fig:hooks-kr}.

On the other hand, the following proposition shows that hook coefficients only fail to distinguish between the $r=1$ and $r=c$ cases.
\begin{proposition}
    Let $G$ and $H$ be two connected unicyclic graphs such that the star-expansions of $\mathbf{X}_G$ and $\mathbf{X}_H$ have the same hook coefficients. Let $r_G$, $r_H$ denote the number of non-trivial rooted trees of each graph, and let $k_G$, $k_H$ denote the respective number of internal edges. Then, up to interchanging $G$ and $H$, either:
    \begin{enumerate}
        \item $r_G=r_H$ and $k_G=k_H$, or
        \item $r_G=1$, $r_H=c$ where $c$ is the size of the cycle in each graph, and $k_G=k_H+1$.
    \end{enumerate}
\end{proposition}
\begin{proof}
    First note that by \cref{cor:cyc size from coeff} both graphs have the same cycle size $c$. Furthermore, the cycle $C_n$ is distinguished among unicyclic graphs by the coefficient $c_{(n)}$ as the unique graph with $c_{(n)}=n-1$. Henceforth, we can assume $r_G, r_H\geq 1$. It follows from \cref{prop:k-and-r-hooks} that $k_G=m+2$ or $k_G=m+1$, where $m$ is the number of 1's in the longest hook appearing in the star-expansion of $\mathbf{X}_G$. We know that if $k_G=k_H$ then $r_G=r_H$ by \cref{prop:k-and-r-hooks}, which takes care of the first case.
    Supposing $k_H\neq k_G$, since $\mathbf{X}_H$, has the same hook coefficients it must be the case that $k_H=m+2$ or $k_H=m+1$ for the same $m$. Thus, without loss of generality let $k_G=m+2$ and $k_H=m+1$. This implies that $r_G=1$ and $r_H \geq 2 $. Moreover, applying the formula of \cref{thm:hook_coeff_general} to the partition $(n-1,1)$ we have
    \begin{align*}
        c_{(n-1,1)}(G)=-[(r_G-1)+(c-1)(k_G-2)]&=-[(r_H-1)+(c-1)(k_H-2)]=c_{(n-1,1)}(H)\\
        r_G+(c-1)m&=r_H+(c-1)(m-1)
    \end{align*}
    so that $r_H-r_G=c-1$, implying $r_H=c$ as desired.
\end{proof}
As a consequence, there is an algorithm to obtain $r$ from the hook coefficients unless $r=1$ or $c$.
\begin{corollary}\label{cor.number.of.rooted.trees}
 We can recover $r$, the number of non-trivial rooted trees attached to the cycle, from the CSF of a connected unicyclic graph whenever $1< r< c$.
\end{corollary}
   \begin{proof}
From \cref{prop:k-and-r-hooks}, we know that the number of internal edges is either $m+1$ or $m+2$, where $m$ is the number of 1's appearing in the longest hook in the support of the star-expansion of $\mathbf{X}_G$. As in the proof of that proposition, substituting $m+1$ and $m+2$ for $k$ into
\[c_{(n-1,1)}=-[(r-1)+(c-1)(k-2)]\]
will give two possible values for $r$ whose difference is $c-1$. Unless these values are $1$ and $c$, one of them will be greater than $c$ (which is impossible) or less than $1$.  If these values are $0$ and $c-1$, whether the graph is $C_n$ ($r=0$) can be determined by the coefficient of $c_{(n)}$, as discussed after \cref{cor:cyc size hooks}. Otherwise, one of these two values will be negative, which is impossible. 
Thus, we only get one valid solution when $1<r<c$ and in that case, we have
\[r=-[(c_{(n-1,1)}-1)+(c-1)(m-1)].\]
\end{proof}

\begin{remark}
As shown in Figure \ref{fig:hooks-kr}, the hook coefficients are not enough to determine the number of nontrivial rooted trees attached to the cycle.  In fact, Corollary \ref{cor.number.of.rooted.trees} is the best we can do, since there are unicyclic graphs with only one rooted tree  attached to a cycle that have the same chromatic symmetric function as another that has $c$ rooted trees attached to the cycle.  Here is the smallest example of an infinite family exhibiting this behavior, as described in \cite{ORELLANA20141}: 
 \begin{figure}[hbt!]
    \centering
    \begin{tikzpicture}[auto=center,every node/.style={circle, fill=black, scale=0.45}, style=thick, scale=0.6] 
        %%%%% Graph 1 %%%%%
    %%%%% Draw vertices %%%%%
    \filldraw[black] (0, 1) coordinate (A0) circle (4pt) node{};
    \filldraw[black] (0, 0) coordinate (A1) circle (4pt) node{};
    \filldraw[black] (1, 0) coordinate (A2) circle (4pt) node{};
    \filldraw[black] (1, 1) coordinate (A3) circle (4pt) node{};
    \filldraw[black] (2, 0) coordinate (A4) circle (4pt) node{};
    \filldraw[black] (2, 1) coordinate (A5) circle (4pt) node{};
    
    %%%%% Draw edges %%%%%
    \draw(A1) -- (A2) -- (A3) -- (A1);
    \draw(A3) -- (A5);
    \draw(A0) -- (A1);
    \draw(A2) -- (A4);

        %%%%% Graph 2 %%%%%
    %%%%% Draw vertices %%%%%
    \filldraw[black, ] (10, 1) coordinate (A0) circle (4pt) node{};
    \filldraw[black] (10, 0) coordinate (A1) circle (4pt) node{};
    \filldraw[black] (11, 0) coordinate (A2) circle (4pt) node{};
    \filldraw[black] (11, 1) coordinate (A3) circle (4pt) node{};
    \filldraw[black] (12, 0) coordinate (A4) circle (4pt) node{};
    \filldraw[black] (12, 1) coordinate (A5) circle (4pt) node{}; %
    
    %%%%% Draw edges %%%%%
    \draw(A0) -- (A1)--(A2)--(A0);
    \draw(A3)--(A5);
    \draw(A2)--(A4);
    \draw(A2)--(A3);
    \end{tikzpicture}
     \caption{Two graphs with the same chromatic symmetric function where one has one nontrivial rooted tree and the other has three.}
     \label{fig:rootedtrees}
\end{figure}

\end{remark}

\section{Paths, cycles, and pans}
\label{sec.special.graphs}

In this section, we explicitly compute the coefficients of all partitions in the star-expansions of paths, cycles, and cycles with a single leaf, also called \defn{pan graphs}. While paths are not unicyclic, they provide a test case for illustration of the technique. All three types of graphs have the advantage that their internal edges can be labeled with consecutive integers $\{1,2,\dots, s\}$ in a natural way. This determines an order in which the star-expansion algorithm is applied, where at each stage we apply the DNC relation to the remaining internal edge with the smallest label. 

In this set-up, we may encode a leaf of the DNC-tree for a graph $G$ as a word $\ul{w}=w_1w_2\dots w_s$ in the alphabet $\{L,M,R,X\}$, where the symbol at position $i$ indicates which branch of the DNC-tree we follow when acting on internal edge $i$. Specifically:
        \begin{itemize}
            \item $w_i=L$ means that when we reach the node of the DNC-tree at which internal edge $i$ is acted upon, we proceed to the graph obtained by deletion of this edge. 
            \item $w_i=M$ means that when acting on internal edge $i$, we proceed to the graph obtained by dot-contraction.
            \item $w_i=R$ means that when acting on internal edge $i$, we proceed to the graph obtained by leaf-contraction. 
            \item $w_i=X$ means that after applying the operations corresponding to $w_1, \dots, w_{i-1}$, the edge with label $i$ is no longer internal so it cannot be acted on. In the DNC-tree, we instead take the next available internal edge and continue the algorithm. 
        \end{itemize}
        
All distinct words that identify valid DNC-tree leaves contribute to the coefficient of some $\stf_\la$ in the star-expansion of $\mathbf{X}_{G}$ in a cancellation-free way. We call a word $\ul{w}$ that produces the star forest $\St_\la$ a \defn{$\la$-word}, and we compute the coefficient $c_\la$ by counting such words. We note that in all of our labelings, internal edges $i$ and $i+1$ will be incident for $1\leq i <s$, which implies that $w_i=L$ will often force $w_{i+1}=X$. On the other hand, the operations $w_i=M$ or $w_i=R$ do not usually impact whether or not $w_{i+1}$ is internal, as we will discuss in the proofs of our results. 

\begin{example}
    Let $G=P_7$ where the internal edges are labeled $\{1,2,3,4\}$ from left to right. Then, $MRML$ and $MLXM$ are two valid $\la$-words for $\la=(3,2,1,1)$, implying that the absolute value of the coefficient $c_\la$ in the star-expansion of $\mathbf{X}_{P_7}$ is at least 2.
\end{example}

\begin{definition}
    If we write a partition as $\la=(\la_1, \dots, \la_{\ell(\la)})$, let $b$ be the index of the last part satisfying $\la_b>1$. We call the partition $\hat{\la}=(\la_1,\dots, \la_b)$ the \defn{body} of $\la$ and $\tilde{\la}=(\la_{b+1}, \dots, \la_{\ell(\la)})=(1,\dots,1)$ the \defn{tail}.
\end{definition}

\begin{proposition} \label{prop:paths}
Let $n\geq 4$ and let $P_n$ the path with $n$ vertices. Then 
\[ \mathbf{X}_{P_n} = \sum_{\lambda= (1^{m_1} 2^{m_2}\ldots n^{m_n})\vdash n} (-1)^{m_1} \binom{m_2+\dots+m_n}{m_2,\dots,m_n} {n-2-\ell(\lambda) + m_1\choose m_1} \mathfrak{st}_\lambda. \]
\end{proposition}
\begin{proof}
Label the internal edges of $P_n$ consecutively from one end to the other by $\{1,\dots, n-3\}$, and apply the star-expansion algorithm to these edges by selecting the edge with smallest possible label at each stage. We analyze each factor that constitutes the coefficient $c_\la$ for $\la=(1^{m_1} 2^{m_2}\ldots n^{m_n})$ separately.
    
Signs in the star-expansion algorithm are introduced only when applying dot-contraction. Starting from a connected graph with at least two vertices, dot-contraction is also the only way to obtain parts of size 1 in a star forest produced by the star-expansion algorithm. Then, $m_1$ must be the number of times $M$ appears in any word that produces $\St_\lambda$, accounting for the factor of $(-1)^{m_1}$. 

We now explain the multinomial coefficient. Observe that any valid $\la$-word $\ul{w}=w_1\dots w_{n-3}$ 
must also satisfy an additional condition: there must be exactly $b-1=\ell(\la)-m_1-1$ symbols $w_i$ with $w_i=L$ in $\ul{w}$, where $b$ is the number of body parts, $\ell(\hat{\la})$. This follows because every deletion of an internal edge produces a new connected component of size at least 2, while neither dot-contraction nor leaf-contraction do so.

The presence of $M$'s anywhere in the string $\ul{w}$ does not affect the part sizes of the body partition so we may ignore them for determining these parts. Hence, let $\ul{v}=v_1\dots v_{n-3-m_1}$ be the word obtained from $\ul{w}$ by omitting all $M$'s. 

Then $\ul{v}$ has the form $\ul{v}=R \cdots R$ if no $L$'s appear, and if at least one $L$ appears then
\[\ul{v}=R\cdots RLXR\cdots RLXR \cdots \cdots R LX R\cdots R,\] where $v_{i_1},\dots ,v_{i_{b-1}}=L$ for some $1\leq i_1 < i_2 < \dots < i_{b-1} \leq n-3-m_1$. In addition, we must have $v_{i_j+1}=X$ if $i_j < n-3-m_1$, since the next edge will no longer be internal. In particular, it may be that $i_{b-1}=n-3-m_1$ in which case this $L$ is not followed by an $X$.

By counting the number of symbols between the $L$'s, we obtain a rearrangement (or ``permutation") $\pi=(\pi_1,\dots, \pi_b)$ of the body parts of $\la$ as follows:
    \[\overbrace{R\cdots R}^{\pi_1-2}\overbrace{LXR\cdots R}^{\pi_2}\overbrace{LXR\cdots R}^{\pi_3}LXR \cdots \cdots R \overbrace{LX R\cdots R}^{\pi_b-1}.\]
That is, $\pi_b = \pi_1 = n-3-m_1$ if there are no $L$'s, and if at least one $L$ appears then
        \begin{itemize}
            \item $\pi_1=i_1+1$,
            \item $\pi_j=i_j-i_{j-1}$ for $2\leq j \leq b-1$ (when $b\geq 3$),
            \item  and $\pi_b=n-1-m_1-i_{b-1}$.
        \end{itemize}

 To account for $\pi_b$, for example, observe that if the last deletion $L$ occurs at internal edge $s$, there is still a path on $n-1-s$ vertices in the connected component containing any remaining internal edges. The number of vertices $\pi_b$ that remain in this last connected component is exactly $n-1-m_1-i_{b-1}$, so that in particular if $i_{b-1}=n-3-m_1$, $\pi_b=2$. Determination of the other parts resulting from the word $\ul{v}$ is similar, noting that each $R$ will grow some part $\pi_j$ by 1 by definition of leaf-contraction.
 
Thus, a valid $\la$-word $\ul{w}$ must satisfy the property that the subword $\ul{v}$ produces a rearrangement $\pi$ of the body parts of $\la$ in this way. Since there is a unique string $\ul{v}$ for each $\pi$, we count $\binom{m_2+\dots+m_n}{m_2,\dots,m_n}$ possibilities.

Once the body part rearrangement is determined, a valid $\la$-word can be completed by inserting $m_1$ symbols $M$ into the word $\ul{v}$ so as to recover a word $\ul{w}$ of length $n-3$.

To count the possibilities here, consider the set of valid DNC-tree words with $b-1$ $L's$ as the disjoint union of two sets.
\begin{enumerate}
    \item $W_1$ consisting of words of length $n-3-(b-1)$ in $\{LX, M, R\}$, where $LX$ counts as one symbol and is used $b-1$ times.
    \item
    $W_2$ consisting of words of length $n-3-(b-1)$ in $\{LX, M, R\}$, where $LX$ is used $b-2$ times, appended by the symbol $L$. 
\end{enumerate}
In either case, we choose which $m_1$ of the $n-3-(b-1)=n-2-(\ell(\la)-m_1)$ symbols should be $M$'s, which can be done in ${n-2-\ell(\lambda) + m_1\choose m_1}$ many ways. Once this is done, the rest of the symbols are determined by the choice of body part rearrangement and associated word $\ul{v}$, so multiplying the binomial coefficient with the multinomial $\binom{m_2+\dots+m_n}{m_2,\dots,m_n}$ gives the number of DNC leaves with star forest $\St_\la$.
\end{proof}

Now we handle the cycle $C_n$, which is, in some sense, the primordial connected unicyclic graph and already gives some indications of how this class differs from trees. As we discussed in the previous section, the coefficients $c_\la$ for hook partitions $\la$ behave differently than non-hook coefficients.\par

\begin{proposition}
\label{prop: cycle_coeff}
Let $C_n$ be the cycle on $n$ vertices for $n\geq 3$, and let 
$\mathbf{X}_{C_n}=\sum_{\la \vdash n} c_\la \mf{st}_\la$. Then for $\lambda = (1^{m_1} 2^{m_2} \dots n^{m_n})$ we have 
    \[ c_\lambda =
    \begin{cases} 
    (-1)^{m_1} \left[ -\binom{n-2}{m_1-1} + (n-1) \binom{n-2}{m_1} \right]
    & \text{if } \la = (n-m_1, 1^{m_1}),\ m_1\neq n-1, 
    \\  \\
   (-1)^{m_1} \frac{n}{\ell(\la)-m_1}\binom{m_2+\dots+m_n}{m_2, \dots ,  m_n} {n-\ell(\lambda)+m_1-1\choose m_1}&\text{else}.\end{cases}\]
   
\end{proposition}

\begin{proof}
    Again label the internal edges $\{1,2,\dots,n\}$ around the cycle and apply DNC relations to the internal edge with smallest label at each stage. As before, each leaf of the DNC-tree is determined by a word in $\{L,M,R,X\}$ of length $n$ which we denote $\ul{w}=w_1\dots w_{n}$. The sign $(-1)^{m_1}$ is always given by the number of $M$'s in a $\la$-word, since both isolated vertices in $\St_\la$ and signs arise only from dot-contraction. We argue for the remainder of each expression separately.\par
\textbf{Hook case:}
The desired equation immediately follows from \Cref{thm:hook_coeff_general}, as $c=n=k$ and $r=0$ in this case. 

\textbf{Non-hook case:} As above, let $\ul{v}=v_1\dots v_{n-{m_1}}$ be the word obtained from a $\lambda$-word $\ul{w}$ by omitting all $M$'s. We will count the contributions from $v_1=L$ and $v_1\neq L$ separately. 

For each body part rearrangement $\pi$ of $\hat{\la}=(\la_1,\dots, \la_b)$, there is a unique corresponding subword $\ul{v}$ yielding $\pi$ for which $v_1=L$. That is, no $R$ appears before the first $L$ so this word has the form
 \[\overbrace{LXR\cdots R}^{\pi_1}\overbrace{LXR\cdots R}^{\pi_2}LXR \cdots \cdots R \overbrace{LX R\cdots RX}^{\pi_b}.\]
To reconstruct a valid $\la$-word $\ul{w}$ with this $\ul{v}$, we choose where to place $m_1$ $M$'s among $n-b-1$ slots since the last symbol is $w_n=X$ here. Thus, these words contribute
    \begin{equation} \label{eq:lxrlx term}
        \binom{m_2+\dots+m_n}{m_2, \dots ,  m_n} {n-\ell(\lambda)+m_1-1\choose m_1}
    \end{equation}
    leaves to $c_\la$. 

    Now we consider $\la$-words with $v_1\neq L$; that is, those $\ul{w}$ for which the subword $\ul{v}$ contains an $R$ before the first $L$. Here, we are free to choose where to place the $M$'s from any of the $n-b$ slots corresponding to letters in the alphabet $\{LX,M, R\}$. 
    The subword $\ul{v}$ can be constructed as
      \[\overbrace{R\cdots R}^{p}\overbrace{LXR\cdots R}^{\pi_1}\overbrace{LXR\cdots R}^{\pi_2}LXR \cdots \cdots R \overbrace{LX R\cdots R}^{\pi_{b-1}}\overbrace{LXR\cdots R}^{\pi_b-p},\]
    so that for any body part rearrangement ending in a given $\pi_b=j$, there are $j-1$ possible subwords yielding the same $\pi$. Thus, this gives a contribution of
    \[\left(\sum_{j=2}^n(j-1) \binom{m_2+\dots+m_n-1}{m_2, \dots , m_j-1,\dots, m_n}\right) {n-\ell(\lambda)+m_1\choose m_1}. \]
    Since 
    \begin{align*}
        \sum_{j=2}^n(j-1) \binom{m_2+\dots+m_n-1}{m_2, \dots , m_j-1,\dots, m_n} 
         &= \sum_{j=2}^n \frac{(j-1)m_j}{m_2+\dots +m_n} \binom{m_2+\dots+m_n}{m_2, \dots, m_n} 
         \\
           &= \frac{1}{m_2+\dots +m_n}\left(\sum_{j=2}^n (j-1)m_j\right) \binom{m_2+\dots+m_n}{m_2, \dots ,  m_n} 
           \\
              &= \frac{n-\ell(\la)}{\ell(\la)-m_1}\binom{m_2+\dots+m_n}{m_2,\dots, m_n}.
 \end{align*}
    Adding with \eqref{eq:lxrlx term} and simplifying gives the result.
\end{proof}

\begin{example} 
By \Cref{prop: cycle_coeff}, the CSF of the $5$-cycle $C_5$ is \[4 \, \mf{st}_{(5)}-11 \, \mf{st}_{(4,1)}+5\, \mf{st}_{(3,2)}+9\, \mf{st}_{(3,1,1)}-5\, \mf{st}_{(2,2,1)}-1\, \mf{st}_{(2,1,1,1)}.\]

For instance, the hook partition $(4,1)$ has $m_1=1$, so $c_{(4,1)}$ is
\[(-1) \left[ -\binom{3}{0} + 4 \binom{3}{1} \right]= -11.\]
For the non-hook partition $(3,2)$, we have $m_1=0$, $m_2=1$, $m_3=1$, and $\ell(3,2)=2$, so $c_{(3,2)}$ is
\[\frac{5}{2}\binom{2}{1,1} {2\choose 0 }=5.\]
Similarly, for $(2,2,1)$ we have $m_1=1$, $m_2=2$, and $\ell(2,2,1)=3$ , so $c_{(2,2,1)}$ is
\[(-1)\frac{5}{2}\binom{2}{2} {2\choose 1}=-5.\]

\end{example}

We highlight a few features of the formula for the star-expansion of $\mathbf{X}_{C_n}$. 
If $m_1=0$ and $\ell(\la)\geq 2$, we can write 
\[c_\lambda = \frac{n(\ell(\lambda) -1)!}{\prod_{i=2}^n m_i!}.\]
This shows that $c_\la=n$ for any two-part partition with $\la_1>\la_2\geq 2$, or $c_{(n/2, n/2)}= n/2$ if $n$ is even. We also see that the leading partition of $\mathbf{X}_{C_n}$ is $(2,1^{n-2})$, since the only smaller partition is $(1^n)$ which is seen to have coefficient 0 by our formula. 
The coefficient of $(2,1^{n-2})$ is always $(-1)^{n-2}$.
\begin{remark}\label{rm.dyckpath}
When $n\geq 6$, for $\mathbf{X}_{C_n}$ our formula gives $c_{(3,3,1^{n-6})}=\binom{n-1}{3}\frac{n+2}{3}$. This sequence begins $0,3,14,40,90,175,\dots$ (OEIS \href{https://oeis.org/A117662}{A117662}) and is equal to the number of certain pairs of non-crossing Dyck paths, or also exterior intersections of diagonals of an irregular $n$-gon. 
\end{remark}

Next we look at the CSF for the unicyclic graph consisting of one $(n-1)$-cycle and one leaf incident to one of the vertices of the cycle, also called the ``pan" on $n$ vertices. 

\begin{example}
The pan graph for $n=6$ and its chromatic symmetric function:
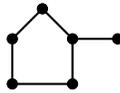
\begin{figure}[hbt!]
    \centering
    \begin{tikzpicture}[auto=center,every node/.style={circle, fill=black, scale=0.45}, style=thick, scale=0.4] 
        %%%%% TREE 1 %%%%%
    %%%%% Draw vertices %%%%%
    \filldraw[black] (0, 0) coordinate (A0) circle (4pt) node{};
    \filldraw[black] (0, 1.5) coordinate (A1) circle (4pt) node{};
    \filldraw[black] (1, 2.5) coordinate (A2) circle (4pt) node{};
    \filldraw[black] (2,1.5) coordinate (A3) circle (4pt) node{};
    \filldraw[black] (2, 0) coordinate (A4) circle (4pt) node{};
    \filldraw[black] (3.5, 1.5) coordinate (A5) circle (4pt) node{};
 
    %%%%% Draw edges %%%%%
   \draw(A0) -- (A1) -- (A2) -- (A3) -- (A4) -- (A0);
    \draw(A3) -- (A5);
\end{tikzpicture}
    \caption{A pan graph}
    \label{fig:pan6}
\end{figure}
 \[4 \,\mf{st}_{(6)}-12\,\mf{st}_{(5,1)}+4\,\mf{st}_{(4,2)}+12\,\mf{st}_{(4,1,1)}+2\,\mf{st}_{(3,3)}-9\,\mf{st}_{(3,2,1)}-4\,\mf{st}_{(3,1,1,1)}+\,\mf{st}_{(2,2,2)}+3\,\mf{st}_{(2,2,1,1)}.\]
\end{example}

\begin{proposition}
        If $G$ is the pan graph on $n$ vertices and $\mathbf{X}_G = \sum_\lambda c_\lambda \mathfrak{st}_\lambda$, 
    then for $\lambda = (1^{m_1} 2^{m_2} \dots n^{m_n})$ we have 
    \[ c_\lambda = \begin{cases} (-1)^{m_1}(n-2){n-3 \choose m_1}, & \text{if } \lambda = (n-m_1, 1^{m_1}),\quad m_1\neq n-1,\\  \\
     (-1)^{m_1}\frac{n-\ell(\lambda)}{\ell(\lambda)-m_1}\binom{m_2+\dots +m_n}{m_2,\dots, m_n} {n-\ell(\lambda)+m_1-1\choose m_1} &\text{else}.\end{cases}\]
    
\end{proposition}
\begin{proof}
    Orient the cycle so that we can label its internal edges from $1$ to $n-1$ starting at one edge that is adjacent to the leaf-edge and ending at the other such edge. We again proceed to apply DNC relations to each internal edge in increasing order of the labels. Denote the word in $\{L,M,R,X\}$ that determines a DNC leaf by $\ul{w}=w_1\dots w_{n-1}$. The sign is given by the number of $M$'s in $\ul{w}$ as in the path and cycle cases.

   \textbf{Hook case:} 
   We immediately obtain the desired equation by \Cref{thm:hook_coeff_general} as $r=1$ and $c=k=n-1$ in this case.
   
    \textbf{Non-hook case:} Here, any $\la$-word must contain at least two $L$'s. We again consider the subword $\ul{v}$ obtained by deleting all $M$'s from a valid $\la$-word. Such $\ul{v}$ then has the form
     \begin{equation} \label{eq:pi pan}
         \overbrace{R\cdots R}^{p}\overbrace{LXR\cdots R}^{\pi_1}\overbrace{LXR\cdots R}^{\pi_2}LXR \cdots \cdots R \overbrace{LX R\cdots R}^{\pi_{b-1}}\overbrace{LXR\cdots R}^{\pi_b-p-1}.
     \end{equation}
    Let $e$ be the leaf-edge of the initial pan graph. In contrast with the cycle case, having no $R$'s before the first $L$ no longer forces $w_{n-1}=X$ because of the presence of $e$. Thus, we can consider all words $\ul{v}$ of the form described in \eqref{eq:pi pan} together, and each one determines a rearrangement $\pi$ of $\hat{\la}=(\la_1,\dots, \la_b)$ as indicated. 

    Notice that there are $\pi_b-1$ ways of constructing a subword $\ul{v}$ that gives the same $\pi$ by just cyclically shifting $\ul{v}$ to vary $p$, the number of initial $R$'s appearing. Thus, for a body part permutation ending with $\pi_b=j$, we have $j-1$ choices of corresponding subword $\ul{v}$. This yields
    \[\sum_{j=2}^n(j-1) \binom{m_2+\dots+m_n-1}{m_2, \dots , m_j-1,\dots m_n} = \frac{n-\ell(\la)}{\ell(\la)-m_1}\binom{m_2+\dots+m_n}{m_2,\dots, m_n}\]
    total subword choices for each body partition $\dot \la=(\la_1,\dots,\la_b)$.
    Once the subword $\ul{v}$ is chosen, we can insert $m_1$ $M$'s by choosing their positions from among $n-\ell(\la)+m_1-1$ slots, and completing the rest with the symbols from $\ul{v}$ to obtain a word of length $n-b-1$ in the alphabet $\{LX,M,R\}$. This accounts for the binomial coefficient in the formula, and multiplying with the number of options for $\ul{v}$ finishes the proof.
\end{proof}

We expect that the techniques of this section can be extended to other families of connected unicyclic graphs. For example, a \defn{tadpole} is a unicyclic graph obtained by attaching a path to a single vertex on the cycle, so the pan graph is a special case where the attached path has length $1$. A \defn{squid} is a connected unicyclic graph with exactly one vertex of degree greater than 2. Alternatively one can think of a squid as a unicyclic graph with one non-trivial rooted tree which consists of some number of paths (called \defn{tentacles}) emanating from the root. These families still have a natural ordering of their internal edges along the cycle and attached paths, so the corresponding DNC-tree leaves can, in principle, be encoded by words similar to the $\{L,M,R,X\}$ words used above. However, for general unicyclic graphs, this technique may not immediately extend. In general, the rooted trees attached to the cycle may branch in many places, there is no single natural ordering of the internal edges, and the dependencies within the DNC-tree become more difficult to control.

% We mention that the techniques of this section could, in principle, be used to construct explicit formulas for star-expansions of other classes of unicyclic graphs like \emph{tadpoles} and \emph{squids}.

\section{Leading partitions for unicyclic graphs}
\label{sec.lead.partition}
In this section, we describe the leading partition of the chromatic symmetric function of any connected unicyclic graph $G$ on $n$ vertices, generalizing a result of Gonzalez, Orellana, and Tomba in~\cite[Theorem~4.15]{gonzalez2024chromatic}. In the case where there are no non-trivial rooted trees attached to the cycle (i.e., $r=0$), the leading partition is $(2,1^{n-2})$ as discussed before \Cref{rm.dyckpath}. In this section we consider the case where $r\geq 1$. We first introduce some notation which we will use throughout this section.

Recall $\sort()$ is the function that takes in a sequence of non-negative integers and outputs it as a partition.
Recall that $I(G)$ is the set of internal edges of the graph $G$ and the connected components of $G \setminus I(G)$ are called the \defn{leaf components} of $G$. 

\begin{definition}
    Notice that $G\setminus I(G)$ is a star forest since every remaining edge is incident with a leaf. We define \defn{$\lc(G)$} to be the partition whose parts are the orders of the connected components of the resulting star forest in $G \setminus I(G)$.
\end{definition}

Recall that $\lead$ was defined in \Cref{def:lead} as the smallest partition, $\lambda$, in lexicographic order such that $c_\lambda \neq 0$ in $\mathbf{X}_G$. 

\begin{theorem}[{\cite[Theorem 4.15]{gonzalez2024chromatic}}]\label{thm:leading.partition.forest}
Let $F$ be a forest. Then $\lead(\mathbf{X}_F) = \lc(F)$.
\end{theorem}

\begin{example}
Let $F$ be the forest in \Cref{{fig:tree.lc}}. After removing the only internal edge $e$ indicated in red, we get $\lead(\mathbf{X}_F) = \lc(F) = (3,2,1)$.
\end{example}

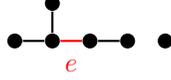
\begin{figure}[H]
\centering
\begin{tikzpicture}[auto=center,every node/.style={circle, fill=black, scale=0.6}, style=thick, scale=0.5] 
\node (V1) at (3,1) {};
\node (V2) at (5,1) {};
\node (V3) at (4,1) {};
\node (t1) at (2,1) {};
\node (t2) at (2,2) {};
\node (t3) at (1,1) {};
\draw (t2) -- (t1);
\draw[color = red] (t1) -- (V1) node[midway, below = 3pt,fill=white] {\huge{$e$}};
\draw (V3) -- (V1);
\draw (t3) -- (t1);
\end{tikzpicture}
\caption{A forest on 6 vertices.}
\label{fig:tree.lc}
\end{figure}

\begin{notation}\label{def.notation}
Let $G$ be a connected $c$-unicyclic graph of order $n$. We label vertices in the $c$-cycle with $v_1,v_2,\dots,v_c$ such that $v_iv_{i+1}$ is an edge for $1\leq i \leq c-1$, as well as $v_cv_1$. Furthermore, each $v_i$ is the root of some rooted tree $T_i$ (which could be trivial).  
Let $\lambda = (\lambda_1, \lambda_2, \dots, \lambda_c)$ be the sequence of positive integers where each part $\lambda_i$ equals $1$ plus the number of leaves attached to $v_i$ in $G$ for $1\leq i \leq c$. Let $\mu^{(i)}$ be the partition obtained from the leaf components of $T_i$ after removing $v_i$ and the leaves incident to $v_i$. Let $\mu = (\mu^{(1)}\cdot \mu^{(2)}\dots \cdot\mu^{(c)})$. 
\end{notation}

Note that for a connected $c$-unicyclic graph $G$ we have $\sort(\lambda \cdot \mu) = \lc(G)$. We will see that this will give the leading partition for $\mathbf{X}_G$ when there are at least two non-trivial rooted trees. When there is only one non-trivial rooted tree, we need to make a minor change to this partition to get the leading partition for $\mathbf{X}_G$.

\begin{figure}[hbt!]
\centering
\begin{tikzpicture}[auto=center,every node/.style={circle, fill=black, scale=0.6}, style=thick, scale=0.5] 
\node[label={[shift={(0.3,-0.8)}]$v_4$}] (V4) at (0,0) {};
\node[label={[shift={(-0.3,-0.8)}]$v_1$}] (V1) at (3,0) {};
\node[label={[shift={(-0.3,-0.2)}]$v_2$}] (V2) at (3,-3) {};
\node[label={[shift={(0.3,-0.2)}]$v_3$}] (V3) at (0,-3) {};
\node (t1) at (-1,0) {};
\node (t2) at (-2,0) {};
\node (t3) at (3,1) {};
\node (t4) at (3,2) {};
\node (t5) at (2,1) {};
\node (t6) at (4,-3) {};
\node (t7) at (3,-4) {};
\node (t8) at (3,-5) {};
\node (t9) at (0,1) {};
\node (t10) at (-1,-3) {};
\node (t11) at (-2,-3) {};
\node (t12) at (4,1) {};
\node (t13) at (4,0) {};
\node (t14) at (-3,-3) {};
\node (t15) at (-2, -4) {};

\draw(t12) -- (V1) -- (V2) -- (V3) -- (V4) -- (V1);
\draw(t9) -- (V4) -- (t1) -- (t2);
\draw(t13) -- (V1) -- (t3) -- (t4);
\draw(t5) -- (t3);
\draw(t8) -- (t7) -- (V2) -- (t6);
\draw(t14) -- (t11) -- (t10) -- (V3);
\draw(t15) -- (t11);
\end{tikzpicture}
\caption{A graph with a 4-cycle.}
\label{eg.4_cycle}
\end{figure}
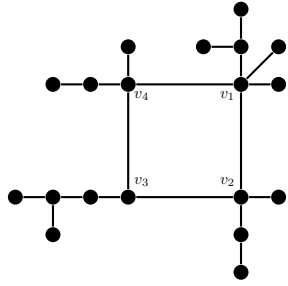

\begin{example}\label{eg.square}
Consider the graph in \Cref{eg.4_cycle}. We have $\lambda = (3,2,1,2)$, $\mu^{(1)} = (3)$, $\mu^{(2)} =(2)$, $\mu^{(3)} = (3, 1)$ and $\mu^{(4)} = (2)$. Hence, $\mu = (3,2,3,1,2)$ and $\lc(G) = \sort(\lambda \cdot \mu) = (3,3,3,2,2,2,2,1,1)$ .
\end{example}

Our main goal for this section is to prove the following.
\begin{theorem} \label{thm:leading}
Let $G$ be a unicyclic graph and $r\geq 1$ denote the number of non-trivial rooted trees attached to the cycle in $G$. Let $\lead = \lead({\bf X}_G)$.

\begin{enumerate}
   \item \label{thm.leading.single} If $r=1$ then  $$\lead=\lc(G \setminus e)$$ where $e$ is chosen to be an edge in the cycle incident to the root of the non-trivial tree.
    \item \label{thm.leading.multiple} If $r\geq 2$ then $$\lead=\lc(G).$$
\end{enumerate} 
\end{theorem}

Note that due to \Cref{lem.disjoint.lead}, it suffices to prove \Cref{thm:leading} for connected unicyclic graphs.

\subsection{One non-trivial rooted tree ($r=1$)}

We will prove the $r = 1$ case of \Cref{thm:leading}(\ref{thm.leading.single}) by induction on the size of the cycle. First we consider the base case, when $G$ contains a triangle.

\begin{lemma}\label{lem.triangle.one.tree}
When $c=3$, \Cref{thm:leading}(\ref{thm.leading.single}) holds.
\end{lemma}
\begin{proof}
Define $\lambda_1,\lambda_2,\lambda_{3}, \mu$ as in \Cref{def.notation}.
First, we note that if $G$ has only one non-trivial rooted tree, say $T_1$ with root $v_1$, $\lc(G \setminus e) = \sort(\lambda_1,\mu,2)$ where $e = v_1v_2$.

The proof follows from a straightforward computation using the DNC relation on $e$ and \Cref{thm:leading.partition.forest}. As $G\setminus e, G \odot e$ and $(G \odot e)\setminus {\ell_e}$ are all forests, \Cref{thm:leading.partition.forest} applies.

\begin{align*}
\lead({\bf X}_{G\setminus e}) & = \lc(G \setminus e) = \sort(\lambda_1,\mu,2)\\
\lead({\bf X}_{G \odot e}) & = \sort(\lambda_1+2,\mu)\\
\lead({\bf X}_{(G \odot e)\setminus {\ell_e}}) & = \sort(\lambda_1+1,\mu,1).
\end{align*}

Since $\sort(\lambda_1,\mu,2) \leq \sort(\lambda_1+1,\mu,1) < \sort(\lambda_1+2 ,\mu)$, we have completed the proof.
\end{proof}

\begin{example}\label{eg.triangle.single}
Let $G$ be the graph in \Cref{eg.3_cycle_1_tree}. We have $\lambda_1 = \lambda_2 = \lambda_3 = 1$ and $\mu = (3)$. Let $e = v_1v_2$. We compute:
\[
\lead({\bf X}_{G\setminus e}) = (3,2,1)\,,\  \lead({\bf X}_{G \odot e})= (3,3)\,,\  \lead({\bf X}_{(G \odot e)\setminus {\ell_e}}) = (3,2,1)\,.
\]
Therefore $\lead({\bf X}_G) = (3,2,1) = \lc(G\setminus e) $.
\begin{figure}[hbt!]
\centering
\begin{tikzpicture}[auto=center,every node/.style={circle, fill=black, scale=0.6},edge/.style = {-,-Latex}, style=thick, scale=0.5] 
\node[label={[shift={(0,-1)}]$v_1$}] (V1) at (2,0) {};
\node[label={[shift={(-0.5,-0.3)}]$v_2$}] (V2) at (4,-3) {};
\node[label={[shift={(0.5,-0.3)}]$v_3$}] (V3) at (0,-3) {};
\node (t1) at (2,1) {};
\node (t2) at (2,2) {};
\node (t3) at (1,1) {};
\draw (V1) -- (V2);
\node [fill = white] at (3.35, -1.2) {\huge $e$};  
\draw (t2) -- (t1) -- (V1);
\draw (V2) -- (V3) -- (V1);
\draw (t3) -- (t1);
\end{tikzpicture}\qquad
\begin{tikzpicture}[auto=center,every node/.style={circle, fill=black, scale=0.6}, style=thick, scale=0.5] 
\node (V1) at (2,0) {};
\node (V2) at (4,-3) {};
\node (V3) at (0,-3) {};
\node (t1) at (2,1) {};
\node (t2) at (2,2) {};
\node (t3) at (1,1) {};
\draw (t2) -- (t1) -- (V1);
\draw (V2) -- (V3) -- (V1);
\draw (t3) -- (t1);
\end{tikzpicture}\qquad
\begin{tikzpicture}[auto=center,every node/.style={circle, fill=black, scale=0.6}, style=thick, scale=0.5] 
\node (V1) at (2,0) {};
\node (V2) at (3,-3) {};
\node (V3) at (0,-3) {};
\node (t1) at (2,1) {};
\node (t2) at (2,2) {};
\node (t3) at (1,1) {};
\draw (t2) -- (t1) -- (V1);
\draw (V3) -- (V1);
\draw (t3) -- (t1);
\draw (V1) -- (V2) node[midway,right = 1pt,fill=white] {\huge{$\ell_e$}};
\end{tikzpicture}\qquad
\begin{tikzpicture}[auto=center,every node/.style={circle, fill=black, scale=0.6}, style=thick, scale=0.5] 
\node (V1) at (2,0) {};
\node (V2) at (3,-3) {};
\node (V3) at (0,-3) {};
\node (t1) at (2,1) {};
\node (t2) at (2,2) {};
\node (t3) at (1,1) {};
\draw (t2) -- (t1) -- (V1);
\draw (V3) -- (V1);
\draw (t3) -- (t1);
\end{tikzpicture}
\caption{Left to right: $G,\ G\setminus e,\  G \odot e,\ (G \odot e)\setminus {\ell_e}$.}
\label{eg.3_cycle_1_tree}
\end{figure}
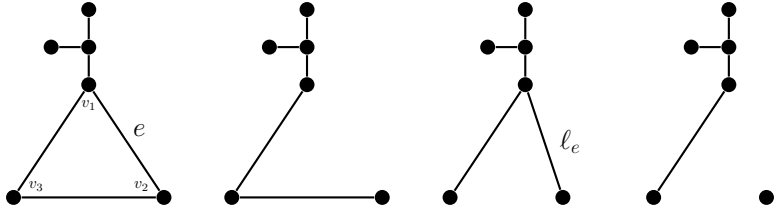
\end{example}

\begin{proposition}
\Cref{thm:leading}(\ref{thm.leading.single}) is true for all $c\geq 3$.
\end{proposition}
\begin{proof}
Define $\lambda, \mu$ as in \Cref{def.notation}.
We prove the statement by induction on $c$. The base case is proven in \Cref{lem.triangle.one.tree}. Suppose the statement holds for any $c$-unicyclic graph with one non-trivial rooted tree where $c \geq 3$. Let $G$ be a $(c+1)$-unicyclic graph with one non-trivial rooted tree, say $T_1$, so that $\lambda_2 = \lambda_3 = \dots = \lambda_{c+1}= 1$. Let $e=v_1v_2$. 

First we compute $\lc(G \setminus e) = \sort(\lambda_1,\lambda_2+\lambda_3,\dots,\lambda_{c+1},\mu) = \sort(\lambda_1,2,1^{c-2},
\mu)$.
Consider the three graphs that we will obtain by applying the DNC relation on the edge $e$ in $G$.
By \Cref{thm:leading.partition.forest}, we have $\lead({\bf X}_{G\setminus e}) = \sort(\lambda_1,2,1^{c-2},
\mu) = \lc(G \setminus e)$.

Since $G \odot e$ is $c$-unicyclic with exactly one rooted tree, we can apply the inductive hypothesis to compute $\lead({\bf X}_{G \odot e})$. Using $e' = v_1v_3$, we get 
\[
\lead({\bf X}_{G \odot e}) = \lc((G \odot e)\setminus e') = \sort(\lambda_1+1,2,1^{c-3},\mu),
\]
as a new leaf becomes incident to $v_1$ as a result of the leaf-contraction.
Similarly, we can compute 
\[
\lead(\mathbf{X}_{(G \odot e)\setminus {\ell_e}}) = \lc(((G \odot e)\setminus {\ell_e})\setminus e')= \sort(\lambda_1,2,1^{c-2},\mu),
\]
where the additional $1$ arises from the dot-contraction.

Since $\lead(\mathbf{X}_{G\setminus e}) = \lc(G\setminus e) = \sort(\lambda_1,2,1^{c-2}, \mu) < \sort(\lambda_1+1,2,1^{c-3},\mu) = \lead({\bf X}_{G \odot e})$, we have completed the proof.
\end{proof}

The following corollaries are evident from the proof given above.

\begin{corollary}\label{cor.not.leaf.contraction}
When we apply DNC relation on a graph satisfying the condition in \Cref{thm:leading}(\ref{thm.leading.single}), then the leading partition does not come from leaf-contraction when the edge is incident to a non-trivial rooted tree.
\end{corollary}

\begin{corollary}
When $G$ is a graph satisfying the condition in \Cref{thm:leading}(\ref{thm.leading.single}), then $2$ must be a part in $\lead(\mathbf{X}_G)$.
\end{corollary}

\subsection{More than one non-trivial rooted tree ($r\geq 2$)}
Here we assume we have more than one non-trivial tree attached at distinct cycle vertices.

\begin{lemma}\label{lem.triangle.multiple}
When $c=3$, \Cref{thm:leading}(\ref{thm.leading.multiple}) holds.
\end{lemma}
\begin{proof}
Since $r\geq 2$, without loss of generality we may assume $T_1$ and $T_2$ are non-trivial trees, and let $e = v_1v_2$.

First note $\lc(G) = \sort(\lambda_1,\lambda_2,\lambda_3,\mu)$. Since in this case $e$ is adjacent to both $T_1$ and $T_2$, this is also equal to $\lead(\mathbf{X}_{G\setminus e})$. We apply the DNC relation on the edge $e$ of $G$.

Next we compute the leading partitions in the other two terms in the DNC relation. Notice that after leaf-contracting, the vertices of the contracted edge are identified to a single vertex, say $v_0$, with $\lambda_1+\lambda_2-1$ leaves. In addition, $G \odot e$ is now a tree, making $\lead(\mathbf{X}_{G \odot e}) = \lc(G \odot e)$. If $T_3$ is a trivial tree, then $\lambda_3 = 1$ and will become a leaf incident to $v_0$ during the leaf-contraction and thus $\lead(\mathbf{X}_{G \odot e}) = \sort(\lambda_1+\lambda_2+\lambda_3,\mu)$. If instead, $T_3$ is not a trivial tree, then the edge $v_0v_3$ is internal and $\lambda_3$ will stay as its own leaf component in $\lc(G \odot e)$. Therefore we have:
\[
\lead(\mathbf{X}_{G \odot e})=
\begin{cases}
\sort(\lambda_1+\lambda_2+\lambda_3,\mu) & \text{ if } T_3 \text{ is trivial}, \\
\sort(\lambda_1+\lambda_2,\lambda_3,\mu) & \text{ if } T_3 \text{ is not trivial}.
\end{cases}
\]

Following a similar discussion,  we can compute $\lead(\mathbf{X}_{(G \odot e)\setminus {\ell_e}})$:

\[
\lead(\mathbf{X}_{(G \odot e)\setminus {\ell_e}})=
\begin{cases}
\sort(\lambda_1+\lambda_2+\lambda_3-1,1,\mu) & \text{ if } T_3 \text{ is trivial}, \\
\sort(\lambda_1+\lambda_2-1,1,\lambda_3,\mu) & \text{ if } T_3 \text{ is not trivial}\,.
\end{cases}
\]

In both cases, we have $\lc(G) = \lead(\mathbf{X}_{G\setminus e}) \leq \lead(\mathbf{X}_{(G \odot e)\setminus {\ell_e}}) < \lead(\mathbf{X}_{G \odot e})$, which proves the statement.
\end{proof}

\begin{example}
Let $G$ be the graph in \Cref{eg.3_cycle_2_trees}. We have $\lambda_1 = 1, \lambda_2 = 3, \lambda_3 = 1$ and $\mu = (3,2)$. In addition $T_3$ is a trivial rooted-tree at $v_3$. Let $e = v_1v_2$. We compute:
\[
\lead(\mathbf{X}_{G\setminus e}) = (3,3,2,1,1)\,,\  \lead(\mathbf{X}_{G \odot e})= (5,3,2)\,,\  \lead(\mathbf{X}_{(G \odot e)\setminus {\ell_e}}) = (4,3,2,1)\,.
\]
Therefore, $\lead(\mathbf{X}_G) = (3,3,2,1,1) = \lc(G)$.
\begin{figure}[hbt!]
\centering
\begin{tikzpicture}[auto=center,every node/.style={circle, fill=black, scale=0.6}, style=thick, scale=0.5] 
\node[label={[shift={(0,-1)}]$v_1$}] (V1) at (2,0) {};
\node[label={[shift={(-0.5,-0.3)}]$v_2$}] (V2) at (4,-3) {};
\node[label={[shift={(0.5,-0.3)}]$v_3$}] (V3) at (0,-3) {};
\node (t1) at (2,1) {};
\node (t2) at (2,2) {};
\node (t3) at (1,1) {};
\node (t4) at (5,-3) {};
\node (t5) at (6,-3) {};
\node (t6) at (5,-2) {};
\node (t7) at (5,-4) {};
\draw (t2) -- (t1) -- (V1);
\draw (V1) -- (V2);
\node [fill = white] at (3.35, -1.2) {\Large{$e$}};
\draw (V2) -- (V3) -- (V1);
\draw (t3) -- (t1);
\draw (V2) -- (t4) -- (t5);
\draw (t7) -- (V2) -- (t6);
\end{tikzpicture}\qquad
\begin{tikzpicture}[auto=center,every node/.style={circle, fill=black, scale=0.6}, style=thick, scale=0.5] 
\node (V1) at (2,0) {};
\node (V2) at (4,-3) {};
\node (V3) at (0,-3) {};
\node (t1) at (2,1) {};
\node (t2) at (2,2) {};
\node (t3) at (1,1) {};
\node (t4) at (5,-3) {};
\node (t5) at (6,-3) {};
\node (t6) at (5,-2) {};
\node (t7) at (5,-4) {};
\draw (t2) -- (t1) -- (V1);
\draw (V2) -- (V3) -- (V1);
\draw (t3) -- (t1);
\draw (V2) -- (t4) -- (t5);
\draw (t7) -- (V2) -- (t6);
\end{tikzpicture}\qquad
\begin{tikzpicture}[auto=center,every node/.style={circle, fill=black, scale=0.6}, style=thick, scale=0.5] 
\node (V1) at (2,0) {};
\node (V2) at (3,-3) {};
\node (V3) at (0,-3) {};
\node (t1) at (2,1) {};
\node (t2) at (2,2) {};
\node (t3) at (1,1) {};

\node (t4) at (3,0) {};
\node (t5) at (4,0) {};
\node (t6) at (3,1) {};
\node (t7) at (3,-1) {};

\draw (V1) -- (V2) node[midway,left = 0.5pt,fill=white] {\Large{$\ell_e$}};
\draw (t2) -- (t1) -- (V1);
\draw (V3) -- (V1);
\draw (t3) -- (t1);
\draw (t4) -- (t5);
\draw (t7) -- (V1) -- (t4);
\draw (V1) -- (t6);
\end{tikzpicture}\qquad
\begin{tikzpicture}[auto=center,every node/.style={circle, fill=black, scale=0.6}, style=thick, scale=0.5] 
\node (V1) at (2,0) {};
\node (V2) at (3,-3) {};
\node (V3) at (0,-3) {};
\node (t1) at (2,1) {};
\node (t2) at (2,2) {};
\node (t3) at (1,1) {};

\node (t4) at (3,0) {};
\node (t5) at (4,0) {};
\node (t6) at (3,1) {};
\node (t7) at (3,-1) {};

\draw (t2) -- (t1) -- (V1);
\draw (V3) -- (V1);
\draw (t3) -- (t1);
\draw (t4) -- (t5);
\draw (t7) -- (V1) -- (t4);
\draw (V1) -- (t6);
\end{tikzpicture}
\caption{Left to right: $G,\ G\setminus e, G \odot e,\ (G \odot e)\setminus {\ell_e}$.}
\label{eg.3_cycle_2_trees}
\end{figure}
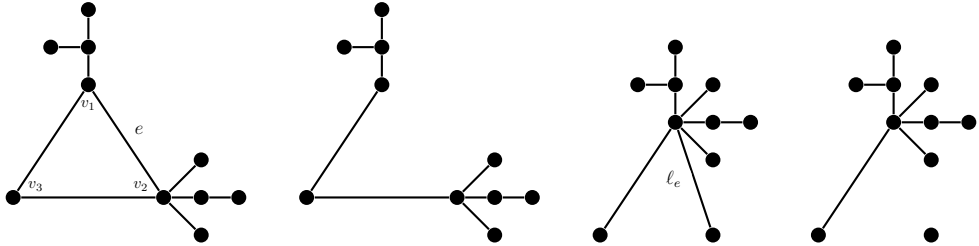
\end{example}

\begin{proposition} \label{prop.thm.leading.rgeq2}
\Cref{thm:leading}(\ref{thm.leading.multiple}) is true for all $c\geq 3$.
\end{proposition}

\begin{proof}
We prove the statement by induction on $c$. The base case is proven in \Cref{lem.triangle.multiple}. Suppose the statement holds for any $c$-unicyclic graph with at least $2$ non-trivial rooted trees for some $c \geq 3$. Let $G$ be a $(c+1)$-unicyclic graph with at least two non-trivial rooted trees. Define $\lambda_1,\dots,\lambda_{c+1},$ and $\mu$ as in \Cref{def.notation}.

We will apply the DNC relation on the edge $e = v_1v_2$, where $T_1$ is a non-trivial tree.
Note that $\lc(G) = \sort(\lambda_1,\dots, \lambda_{c+1},\mu)$.
In addition, we see that in $G \odot e$ (resp. $(G \odot e)\setminus {\ell_e}$), the vertices of the contracted edge are identified to a single vertex, say $v_0$, with $\lambda_1+\lambda_2-1$ (resp. $\lambda_1+\lambda_2-2$) leaves incident to $v_0$. To compute $\lead(\mathbf{X}_{G \odot e})$ and $\lead(\mathbf{X}_{(G \odot e)\setminus {\ell_e}})$, we discuss the following two cases.

\begin{enumerate}
\item When $G \odot e$ is $c$-unicyclic with at least two non-trivial rooted trees, so is $(G \odot e)\setminus {\ell_e}$. In this case, we apply our inductive hypothesis and get
\[\lead(\mathbf{X}_{G \odot e}) = \lc(G \odot e) = \sort(\lambda_1+\lambda_2,\lambda_3,\dots,\lambda_{c+1},\mu)\,, \text{ and}
\]
\[
\lead(\mathbf{X}_{(G \odot e)\setminus {\ell_e}})= \lc((G \odot e)\setminus {\ell_e}) = \sort(\lambda_1+\lambda_2-1,1,\lambda_3,\dots, \lambda_{c+1},\mu).
\]
\item When $G \odot e$ is $c$-unicyclic with exactly one non-trivial rooted tree, so is $(G \odot e)\setminus {\ell_e}$. This happens exactly when $T_1, T_2$ are non-trivial and $T_k$ is trivial for all $3\leq k\leq c+1$. In this case, we have $\lambda_3 = \dots = \lambda_{c+1} = 1$. Let $e' = v_0v_3$. We apply \Cref{thm:leading}(\ref{thm.leading.single}) and get 
\begin{align*}
\lead(\mathbf{X}_{G \odot e}) &= \lc((G \odot e)\setminus e') \\
& = \sort(\lambda_1+\lambda_2,\lambda_3 + \lambda_4,\lambda_5,\dots,\lambda_{c},\lambda_{c+1},\mu) = \sort(\lambda_1+\lambda_2,2, 1^{c-3},\mu)\,, \text{ and}\\
\lead(\mathbf{X}_{(G \odot e)\setminus {\ell_e}}) &= \lc((G \odot e)\setminus {\ell_e})\setminus e')\\
& = \sort(\lambda_1+\lambda_2-1,1,\lambda_3 + \lambda_4,\lambda_5,\dots, \lambda_{c},\lambda_{c+1},\mu) = \sort(\lambda_1+\lambda_2-1,2,1^{c-2},\mu).
\end{align*}
\end{enumerate}
Note that in both cases we have $\lead(\mathbf{X}_{(G \odot e)\setminus {\ell_e}}) < \lead(\mathbf{X}_{G \odot e})$.

Finally we compute $\lead(\mathbf{X}_{G\setminus e})$ and compare it with $\lead(\mathbf{X}_{(G \odot e)\setminus {\ell_e}})$. We discuss the following two cases.
\begin{enumerate}
\item When $T_2$ is a non-trivial rooted tree, then the edge $v_2v_3$ is internal in $G\setminus e$ and we have 
\[
\lead(\mathbf{X}_{G\setminus e}) = \lc(G) = \sort(\lambda_1,\dots, \lambda_{c+1},\mu) \leq \lead(\mathbf{X}_{(G \odot e)\setminus {\ell_e}}).
\]
\item When $T_2$ is a trivial-rooted tree, then the edge $v_2v_3$ is a leaf-edge in $G\setminus e$. In this case, $\lambda_2 = 1$ and $(G \odot e)\setminus {\ell_e}$ is $c$-unicyclic with at least two non-trivial rooted trees.
We have
\begin{align*}
\lead(\mathbf{X}_{G\setminus e}) &= \sort(\lambda_1,\lambda_2+\lambda_3,\lambda_4,\dots,\lambda_{c+1}, \mu) = \sort(\lambda_1,1+\lambda_3,\lambda_4,\dots,\lambda_{c+1}, \mu) \\
&> \sort(\lambda_1,1,\lambda_3,\dots,\lambda_{c+1},\mu) = \lc(G) = \sort(\lambda_1+\lambda_2-1,1,\lambda_3,\dots, \lambda_{c+1},\mu) = \lead(\mathbf{X}_{(G \odot e)\setminus {\ell_e}}).
\end{align*}
\end{enumerate}
Therefore the leading partition of $G$ must be $\lc(G)$.
\end{proof}
The following corollary is evident from the proof of \cref{prop.thm.leading.rgeq2}, and will be used in the next section.
\begin{corollary}\label{cor.leading.comparison}
    When $r\geq 2$ and DNC relation is applied to an edge $e$ which is incident to one trivial tree root and one non-trivial tree root, then $\lead(\mathbf{X}_{(G\odot e)\setminus \ell_e})< \lead(\mathbf{X}_{G \odot e})$ and $\lead(\mathbf{X}_{(G\odot e)\setminus \ell_e}) < \lead(\mathbf{X}_{G\setminus e})$.
\end{corollary}

\begin{example}
Consider the graph $G$ from \Cref{eg.square}. Let $e = v_1v_2$. Recall that  $\lambda = (3,2,1,2)$, $\mu = (3,2,3,1,2)$. We compute that
\begin{align*}
\lc(G) = (3,3,3,2,2,2,2,1,1) &= \lead(\mathbf{X}_{G\setminus e}) < \lead(\mathbf{X}_{(G \odot e)\setminus {\ell_e}}) = (4,3,3,2,2,2,1,1,1) \\
&< \lead(\mathbf{X}_{G \odot e}) = (5,3,3,2,2,2,1,1).
\end{align*}
\end{example}

\begin{corollary}
\label{cor.num.leaves}
Let $\lead$ be the leading partition in the star-expansion of $\mathbf{X}_G$ where $G$ is a connected unicyclic graph with a cycle of size $c$. If $G$ has exactly $1$ non-trivial rooted tree, then the number of leaves in $G$ is equal to $|\lead| - \ell(\lead)-1$; if $G$ has at least $2$ non-trivial rooted trees, then the number of leaves in $G$ is equal to $|\lead| - \ell(\lead)$.
\end{corollary}

\begin{proof}
It follows from the proofs of \Cref{thm:leading}(\ref{thm.leading.single}) and \Cref{thm:leading}(\ref{thm.leading.multiple}) that 

\[
\lead(\mathbf{X}_G)=
\begin{cases}
\sort(\lambda_1,\dots, \lambda_{c-2},\lambda_{c-1}+\lambda_c,\mu) & \text{ if } r = 1,  \\
\sort(\lambda_1,\dots,\lambda_c,\mu) & \text{ if } r>1,
\end{cases}
\]
where we recall $\lambda_i-1$ is the number of leaves attached to $v_i$ on the cycle and $\mu$ is the concatenation of leaf components not containing the cycle vertices.
\end{proof}

\subsection{Application to cuttlefish graphs}
\label{sec.applications}

% A \defn{squid} is a connected unicyclic graph with exactly one vertex of degree greater than 2. Alternatively one can think of a squid as a unicyclic graph with 1 non-trivial rooted tree which consists of some number of paths (called \defn{tentacles}) emanating from the root. 
Recall that a squid is a connected unicyclic graph with exactly one vertex of degree greater than $2$.
We define \defn{cuttlefish} to be squids such that all tentacles are of length 1.
Let $\mc{U}_n$ denote the set of connected unicyclic graphs on $n$ vertices, and let $\mc{U}_{c,n}$ denote the set of such graphs whose cycles are of order $c$.

\begin{proposition}
\label{prop: cuttefish_distinguished}
    The chromatic symmetric function can distinguish cuttlefish from among all connected unicyclic graphs on $n$ vertices. 
\end{proposition}
\begin{proof} Let $C_{c,t}$ denote the cuttlefish on $c+t$ vertices which is constructed by attaching $t$ leaves to a single vertex of a $c$-cycle. Since by Corollary~\ref{cor:cyc size from coeff} we can obtain the size of the cycle of a connected unicyclic graph from the coefficient of $c_{(n)}$, it is sufficient to show that $\mathbf{X}_{G}$ distinguishes $C_{c,t}$ from other graphs in $\mc{U}_{c,c+t}$. Note that $C_n$ is the only graph in $\mc{U}_{c,c}$, so we henceforth suppose $t>0$.

 By \Cref{thm:leading}(\ref{thm.leading.single}), the leading partition of $\mathbf{X}_{C_{c,t}}$ in the star-basis is $\lead=(t+1,2,1^{c-3})$. We will show that $C_{c,t}$ is the only graph in $\mc{U}_{c,n}$ that has this leading partition.

Let $G\in \mc{U}_{c,c+t}$ with $G\not\cong C_{c,t}$. If $r=1$, then by \Cref{thm:leading}(\ref{thm.leading.single}), $\lead(\mathbf{X}_G) = (\lambda_1,2,\mu,1^{c-3})$ where  $\lambda_1-1$ is the number of leaves attached to the root and $\mu$ is the leaf components in the remaining part of the tree. If $\lead(\mathbf{X}_G) = (t+1,2,1^{c-3})$, then $\mu$ is empty and $\lambda_1-1 = t$, making $G = C_{c,t}$. If $r\geq 2$, then each non-trivial tree has at least one leaf component of size at least two. To match the first two entries of $\lead(\mathbf{X}_{C_{c,t}})$, the only possibility is $r=2$ where one of the trees is $\St_{t+1}$, another is $\St_2$, and the rest are trivial. But now there are $c-2$ cycle vertices remaining as their own leaf components, so $|G|=c+t+1$, which is a contradiction. This completes the proof.
\end{proof}
\begin{remark}
    Martin, Morin, and Wagner show in \cite[Theorem 12]{martin2007distinguishingtreeschromaticsymmetric} that no two squids have the same CSF. Note that while the class of squids is larger than that of cuttlefish, their result differs from the proposition above in that they do not distinguish squids from among all unicyclic graphs. On the other hand, both results implicitly give a procedure for reconstructing the graph from the CSF when the ambient class of graphs is known.
\end{remark}
%%%%%%%%%%%%%%%%%%%%%%%%%%%%%%%%%%%%%%%%%%%%%%%%%%%%%%%
\section{Leading coefficients for unicylic graphs}
\label{sec.lead.coeff}
In this section, we provide formulas for the leading coefficient of any unicyclic graph of order $n$. The formulas depend on the number of rooted trees in the unicyclic graph, denoted by $r$. In the case where $r=0$, the leading coefficient is $(-1)^{n-2}$, as discussed following Proposition~\ref{prop: cycle_coeff}. We now recall a result from \cite{gonzalez2024chromatic} on the leading coefficient for trees, which will be useful in our extensions to the unicyclic case. Recall that a vertex is called deep if it is internal and has no leaf as a neighbor.
\begin{theorem}[{\cite[Theorem 4.28]{gonzalez2024chromatic}}]
\label{thm:lead.coeff.tree}
    Suppose that $T$ is a tree of order $t$ with deep vertices $\{u_1,\dots, u_p\}$ such that $d_i=\deg(u_i)$ for all $i$. Then, the leading coefficient of $\mathbf{X}_T$ is \[(-1)^p \prod_{i=1}^p (d_i-1).\]
\end{theorem} 

\Cref{thm:lead.coeff.tree} can be reformulated using the elementary symmetric functions, see \Cref{section:Sym}. 
\begin{remark}
     \emph{Elementary symmetric polynomials} in a finite number of variables are defined analogously to the elementary symmetric functions.  Observe that we have the identity
    \[\prod_{i=1}^p(d_i-1)=\sum_{j=0}^p(-1)^j e_{p-j}(d_1,\dots,d_p).\]
    That is, the leading coefficient is an alternating sum of elementary symmetric polynomials in the degrees of the deep vertices.
\end{remark}

\subsection{One non-trivial rooted tree ($r=1$)}
We first give a formula for the leading coefficient of any unicyclic graph when $r=1$ depending on whether or not the root is a sprout, that is, a deep vertex on the cycle of degree at least 3.
In this subsection, we will use $\varepsilon_i$'s to denote edges since we will use $e_i$'s for the elementary symmetric functions.
\begin{theorem}
    Let $G$ be a connected unicyclic graph with a cycle of size $c$ and exactly one non-trivial rooted tree $T$. Let $v$ be the root of $T$. Label the edges on the cycle within $G$ with edges $\varepsilon_1$ and $\varepsilon_c$ incident to $v$, so that the labels increase $\varepsilon_2, \varepsilon_3, \dots , \varepsilon_c$ in a clockwise direction. Consider the tree $T'$ obtained by deleting edge $\varepsilon_1$. Suppose that $T'$ has $p$ deep vertices $\{u_1,\dots, u_p\}$ such that $d_i=\deg(u_i)$ for all $i$. 
   \begin{enumerate}
       \item If the root $v=u_1$ is a sprout in $G$, then 
       \[c_\lead= (-1)^p\left[(c-2)\prod_{i=1}^p (d_i-1)+\prod_{i=2}^p (d_i-1)\right].\]
       \item If the root $v=u_1$ is not a sprout in $G$, then 
       \[ c_\lead= (-1)^p(c-2)\prod_{i=1}^p (d_i-1).\]
   \end{enumerate} 
\end{theorem}

\begin{proof}
    First, we consider case (1) by identifying graphs from the leaves in the DNC-tree that have the same $\lead$ as $\mathbf{X}_G$. By the fact that the star-expansion algorithm provides a cancellation-free formula, it suffices to determine the contribution to $c_\lead$ from each of these graphs. We can ignore the terms obtained through a leaf-contraction of $\varepsilon_1$ since, by \Cref{cor.not.leaf.contraction}, it will not have the same leading partition as $\mathbf{X}_G$. Upon deleting $\varepsilon_1$, we obtain the tree $T'$ where $\lead(\mathbf{X}_{T'})= \lead(\mathbf{X}_G)$ by \Cref{thm:leading}(\ref{thm.leading.single}). 
    By \Cref{thm:lead.coeff.tree}, the contribution from the deletion of $\varepsilon_1$ to $c_\lead$ is \[(-1)^{p}\prod_{i=1}^p(d_i-1).\]
    Now, consider $G'=(G \odot \varepsilon_1)\setminus {\ell_{\varepsilon_1}}$.
    If $c=3$, then $G'$ is a tree together with an isolated vertex and $v$ is no longer a deep vertex; it has one leaf incident to it. Hence, the cycle's vertices contribute one part of size 2 (leaf component containing $v$) and one part of size 1 to $\lead$ according to \Cref{thm:leading}(\ref{thm.leading.single}), and the other leaf components of the graph are unaffected. Thus, the resulting forest has the same leading partition as $\mathbf{X}_G$, and its contribution to $c_\lead$ is 
    \[-(-1)^{p-1}\prod_{i=2}^p (d_i-1),\] by \cref{thm:lead.coeff.tree}, with the additional minus sign coming from the dot-contraction in the DNC relation applied to $\varepsilon_1$. Combining with the above gives the claimed formula when $c=3$. 

    If $c>3$, $G'$ is a unicyclic graph with one non-trivial rooted tree and an isolated vertex, and by \Cref{thm:leading}(\ref{thm.leading.single}) $\lead(G') = \lc(G'\setminus \varepsilon_c)$. Letting $v'$ be the other vertex on the edge $\varepsilon_1$, then $G'$ is obtained from $G$ by deleting $\varepsilon_1$, identifying $v$ and $v'$ into a single vertex $v_0$ (which is connected to the tree), and adding an isolated vertex $v^*$. We show that $\lead(G') = \lc(G' \setminus \varepsilon_c) = \lc(G \setminus \varepsilon_c) = \lead(G) $. First we note the leaf components for the non-trivial rooted tree (including $v$ in $G\setminus \varepsilon_c$ and $v_0$ in $G'\setminus \varepsilon_c$) are equal in $G'\setminus \varepsilon_c$ and $G\setminus \varepsilon_c$. Second, for the vertices on the cycle (excluding $v$ in $G\setminus \varepsilon_c$ and $v_0$ in $G'\setminus \varepsilon_c$), $G \setminus \varepsilon_c$ has leaf components of orders $(2,1^{c-3})$ and $G'\setminus \varepsilon_c$ has leaf components of orders $(2,1^{c-4})$ from the cycle vertices. Collecting these together with the isolated vertex in $G'$, we see that $\lead(G')=\lead(G)$.
    
    To get $c_\lead$, we must continue down the DNC-tree by next applying the relation to $\varepsilon_2$ in $G'$. Again, graphs arising from leaf-contraction will have CSFs with lexicographically larger leading partitions, so we can ignore them. Deletion of the edge leads to a forest with the same $\lead$ and $c_\lead$ as $\mathbf{X}_{T'}$ but with one fewer deep vertices of degree two. Together with the sign contributed by the previous dot-contraction, this yields another term contributing \[(-1)^{p}\prod_{i=1}^p(d_i-1).\]
    To $G'$, we continue to apply the DNC relation to the edges $\varepsilon_3, \varepsilon_4,\ldots, \varepsilon_{c-2}$ in consecutive order. At each step, we gain $(-1)^{p}\prod_{i=1}^p(d_i-1)$ through deletion. In the last dot-contraction, we gain 
    \[(-1)^{c-2}(-1)^{p-(c-2)}\prod_{i=2}^p(d_i-1).\]  This term is analogous to the $c=3$ case, where $(-1)^{c-2}$ represents the number of dot-contractions and $(-1)^{p-(c-2)}$ corresponds to the remaining deep vertices at the final step. Combining all these DNC-tree contributions gives the claimed formula for case (1).    

    The proof for case (2) follows a similar approach to case (1). Suppose $v$ is not deep in $G$. Then, the deletion of $\varepsilon_1$ produces a tree $T'$ with leading coefficient exactly  \[(-1)^{p}\prod_{i=1}^p(d_i-1).\]
    If $c=3$, both dot-contraction and leaf-contraction of $\varepsilon_1$ increase the size of the leaf component containing $v$ by one, so their leading partitions are different from $\mathbf{X}_{T'}$ and thus do not contribute to $\lead$. On the other hand, if $c>3$, dot-contraction does not increase the size of the leaf component containing $v$, and we may continue applying the DNC relations on $\varepsilon_2, \varepsilon_3, \dots, \varepsilon_{c-2}$. Each deletion contributes the same term as above to $c_\lead$, occurring a total of $(c-2)$ times. The final dot-contraction on edge $\varepsilon_{c-2}$ does not contribute to $c_\lead$, as it will change the leading partition of the CSF of the resulting forest, just as in the $c=3$ case.
    \end{proof}

\begin{example}
Consider $G$ as in \Cref{eg.triangle.single} and \Cref{eg.3_cycle_1_tree}, which has one non-trivial rooted tree and $v_1$ is deep in $G\setminus e$. We get $d_1 = 2, \ p = 1$, and  $c=3$. We compute 
\[
c_{\lead} = (-1)^1 \left[ (3-2) \cdot 1 + 1\right] = -2\,.
\]
This is consistent with 
\[
\textbf{X}_G = 2\mathfrak{st}_{(6)} -4\mathfrak{st}_{(5,1)}+\mathfrak{st}_{(4,2)}+2\mathfrak{st}_{(4,1,1)}+2\mathfrak{st}_{(3,3)}-2\mathfrak{st}_{(3,2,1)}.
\]
\end{example}

\subsection{More than one non-trivial rooted tree ($r\geq 2$)}
Here, the formulas again depend on whether or not the root vertices of non-trivial trees on the cycle are deep or not. 
\begin{definition}
    Let $G$ be a unicyclic graph with more than one non-trivial tree at distinct cycle vertices.  Let the set of sprouts be denoted $v_1,\dots, v_s$, and their degrees $b_i=\deg(v_i)$. Let the non-sprout deep vertices be denoted $u_1,\dots,u_p$ with $d_i=\deg(u_i)$ as before. 
\end{definition}
\begin{proposition} \label{prop:no sprouts clead}
    Let $G$ be a connected unicyclic graph with no sprouts and $r\geq 2$. Then 
    \[c_{\lead}=(-1)^p\prod_{i=1}^p(d_i-1).\]
\end{proposition}

\begin{proof}
    If there is a cycle edge incident to two non-trivial tree roots, then applying DNC to that edge shows that the only contribution to the leading coefficient will come from the tree obtained by deletion. If there is no such edge, then the only contribution to the leading coefficient will come from dot-contraction, so we dot-contract until we have two adjacent non-trivial tree roots and then same argument applies. The sign is the same in both cases.
\end{proof}

Note that unlike the $r=1$ case, deletion of edges that are not incident to two roots of non-trivial trees no longer produces a tree with the same leading partition, which explains the absence of the factor $c$ in the formula of \cref{prop:no sprouts clead}.
If we have at least one sprout, then our formula will depend on both the number of sprouts, which we denote by $s$, as well as $r$, the number of non-trivial trees in total.

\begin{theorem}
Let $b_i$ ($1\leq i \leq s$) be the degrees of sprout vertices, and $d_i$ ($1\leq i \leq p$) be the degrees of non-sprout deep vertices of a connected unicyclic $G$ with $r\geq 2$ and $c\geq 3$. Then
\[c_\lead=(-1)^{p+s}
\begin{cases}
\displaystyle\prod_{i=1}^p(d_i-1)\,\left[\left(\prod_{i=1}^s(b_i-1)\right)-\sum_{i=1}^s b_i+2s-1\right] & \text{ if  } r=s,\\
\\
\displaystyle\prod_{i=1}^p(d_i-1)\,\left[\left(\prod_{i=1}^s(b_i-1)\right)-1\right] & \text{ if  } r=s+1,\\
\\
\displaystyle\prod_{i=1}^p(d_i-1)\, \prod_{i=1}^s(b_i-1) & \text{ if  } r\geq s+2.\\
\end{cases}\]
\end{theorem}

\begin{proof}
   We proceed by induction on $r$.
   \noindent\textit{Base cases ($r=2,3$):}  First suppose that $r=2$. If $c>3$, this means that there is a cycle edge $e$ incident to the root of one trivial tree and one non-trivial tree. Since $r= 2$, by Theorem~\ref{thm:leading} (\ref{thm.leading.multiple}) $\lead(\mathbf{X}_G)=\lc(G)$. Then, by \cref{cor.leading.comparison}, both deletion and leaf-contraction of $e$ will produce graphs whose CSFs will have larger leading partitions, meaning we can ignore these terms for the purpose of computing $c_\lead(\mathbf{X}_G)$ and consider only the CSF of the graph obtained by the dot-contraction of $e$. By \cref{lem.disjoint.lead}, the leading partition of $\mathbf{X}_G$ and $\mathbf{X}_{(G\odot e)\setminus \ell_e}$ are the same, and we may repeat this argument on any cycle edge of $(G\odot e)\setminus \ell_e$, the dot-contracted graph, which is incident to a cycle vertex of degree two and the root of a non-trivial tree.
   Thus we continue to dot-contract edges of $G$ until we obtain a graph which consists of the union of a unicyclic component $G'$ which contains a triangle ($c'=3$) and two non-trivial trees, together with some isolated vertices. 
    The only impact of these dot-contractions on $c_{\lead}$ is the sign, which is accounted for in the factor $(-1)^p$ in the formula above. This is because we will have dot-contracted $c-3$ times to obtain a triangle, and removed the same number of non-sprout deep vertices from $G$ in the process.
    
Let $e'$ be the edge of $G'$ between the roots of the non-trivial trees. The case $s=0$ is handled by \cref{prop:no sprouts clead}. If $s=1$, then both dot-contraction and leaf-contraction of $e'$ give different leading partitions than $\lead(\mathbf{X}_G)$, so we can just consider $c_\lead(\mathbf{X}_{G'\setminus e'})$. This graph is now a tree for which all of the non-sprout deep vertices have the same degree as in $G$, and the unique vertex that was a sprout with degree $b_1$ in $G$ now has degree $b_1-1$. Applying the tree formula of \cref{thm:lead.coeff.tree} to $G'\setminus e'$ recovers the second case of the claimed formula above. Finally if $s=2$, we can again ignore dot-contraction and leaf-contraction, but now deletion of $e'$ gives that the sprout vertices enter into the tree formula as a factor of
\[((b_1-1)-1)((b_2-1)-1)=(b_1-1)(b_2-1)-(b_1+b_2)+3=(b_1-1)(b_2-1)-(b_1+b_2)+2s-1\] 
    as claimed, completing the case $r=2$.
Now suppose $r=3$. As in the $r=2$ case, we may dot-contract until we obtain $G'$ whose connected unicyclic component consists of a triangle with the same three rooted trees attached and $|c_\lead(\mathbf{X}_{G'})|=|c_{\lead}(\mathbf{X}_G)|$. If $s\leq 1$, then we let $e'$ the edge between two non-sprout roots, and the deletion term of DNC applied to $e'$ gives the only contribution to $c_{\lead}$. As $G'\setminus e'$ is a tree, the claimed formula holds.

Next let $s=2$, and let $e'$ be the edge between sprouts $v_1$ and $v_2$ in $G'$. Applying DNC to $e'$, we see that both deletion and dot-contraction terms contribute to $c_\lead$. For brevity, we only discuss the ``sprout factors" in our formula, as the degrees of the other deep vertices of the graph are unchanged by DNC operations on cycle edges. Deletion yields a tree  with sprout factor $(b_1-2)(b_2-2)$,
while dot-contraction merges the sprouts into a single deep vertex of degree $b_1+b_2-3$ since the triangle collapses, yielding another tree. Adding the sprout factor contributions of these trees gives
\[b_1b_2-2(b_1+b_2)+4+(b_1+b_2-4)=(b_1-1)(b_2-1)-1\]
as claimed. 

If $s=3$, suppose we apply DNC to the edge between sprouts $v_1$ and $v_2$. Again we must consider the deletion and dot-contraction terms.  Now, the sprout factor contributions of the corresponding trees combine to give
    \[(b_1-2)(b_2-2)(b_3-1)+(b_1+b_2-4)(b_3-2)=b_1b_2b_3-(b_1b_2+b_1b_3+b_2b_3) +4.\]
    We see that this sprout factor is the same as 
    \[(b_1-1)(b_2-1)(b_3-1)-(b_1+b_2+b_3)+5\]
    as claimed. This completes the base cases.

\noindent\textit{Inductive step:}
    Now suppose $G$ has $r>3$, and that our formula for $c_\lead$ holds for unicyclic graphs with $r-1$ non-trivial trees. As in the $r=3$ case, by \cref{lem.disjoint.lead} when $c>r$ we may apply dot-contraction to edges of the cycle until we obtain a unicyclic graph $G'$ with cycle size $c'=r$ and such that $|c_\lead(\mathbf{X}_{G'})|=|c_{\lead}(\mathbf{X}_G)|$.
    
    We consider three cases for the graph $G'$.
    \begin{enumerate}[label=(\alph*)]
        \item There are two adjacent non-sprout cycle vertices.
        \item There are two adjacent sprouts.
        \item Every pair of neighbors on the cycle consists of sprout and a non-sprout.
    \end{enumerate}\par    
    \textbf{Case (a):} This case is easiest. First notice that this assumption implies $r\geq s+2$. Let $e'$ be an edge of $G'$ between two non-sprouts, and apply DNC to this edge. As above, we only need to consider the contribution of $G'\setminus e'$ which is now a tree in which all deep vertices have the same degree as they do in $G$. Thus, we recover the third case of our claimed formula for $c_{\lead}$.\par
    \textbf{Case (b):} 
 Let $e'$ be an edge between sprouts $v_1$ and $v_2$. The deletion of $e'$ produces a tree with the same leading partition and sprout factor contribution
 \begin{equation}\label{eq:G' minus e'}
 (b_1-2)(b_2-2)\prod_{i=3}^s(b_i-1),
 \end{equation}
 while dot-contraction produces a unicyclic graph $G''$ with number of non-trivial trees $r''=r-1$ and number of sprouts $s''=s-1$. 

 If $r=s$, then $r''=s''$ and the merging of $v_1$ and $v_2$ produces a sprout of degree $b_1+b_2-2$. Hence, by the inductive hypothesis, the sprout factor contribution of $G''$ is 
 \[(b_1+b_2-3)\left(\prod_{i=3}^s(b_i-1)\right)-(b_1+b_2-2+b_3+\dots+b_s)+2s-3.\]
Adding with the contribution of $G'\setminus e'$ gives the desired formula.

If $r=s+1$, the sprout factor contribution of $G'\setminus e'$ is unchanged. On the other hand, dot-contraction gives a sprout factor of
\[(b_1+b_2-3)\left(\prod_{i=3}^s(b_i-1)\right)-1,\]
by the inductive hypothesis. Adding the quantity we obtained in \eqref{eq:G' minus e'} gives the desired result. 

Finally, if $r \geq s+2$ then $G''$ contributes sprout factor
\[(b_1+b_2-3)\prod_{i=3}^s(b_i-1),\]
and again adding with \eqref{eq:G' minus e'} gives the result.

\textbf{Case (c):} Note that the condition of this case together with $r\geq 4$ implies that $r\geq s+2$. Let $e'$ be an edge between sprout $v_1$ and another cycle vertex which is adjacent to $v_1$. 
We again consider sprout factor contributions of the deletion and dot-contraction terms of DNC applied to $e'$. $G'\setminus e'$ has a sprout factor of 
\[(b_1-2)\prod_{i=2}^s(b_i-1)\]
while the dot-contraction gives a graph $G''$ with $r''=r-1$, $s''=s-1$ to which we apply the inductive hypothesis. As $r''\geq s''+2$, the sprout factor is just
\[\prod_{i=2}^s(b_i-1),\]
and adding these terms gives the result. This concludes this case, the inductive step, and the proof of the theorem.
\end{proof}

\begin{example}
Consider $G$ in \Cref{eg.prelim}, which has $r=3$ non-trivial rooted trees. The only sprout is $v_3$ with degree $3$ and the non-sprout deep vertices are $v_4$ and $v_{11}$, both with degree $2$. So $s = 1, p = 2$ and we get
\[
c_{\lead} = (-1)^{2+1}\cdot 1 \cdot 1 \cdot 2 = -2\,.
\]
\end{example}

\section{Bicyclic graphs}
\label{sec.general.graphs}

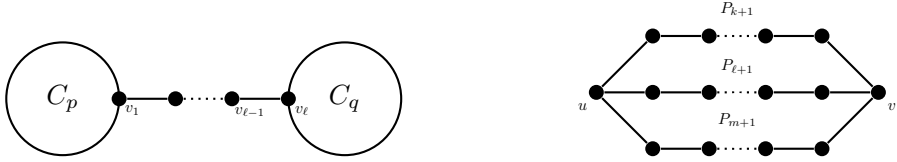
\begin{figure}[hbt!]
\centering
\begin{tikzpicture}[auto=center,every node/.style={circle, fill=black, scale=0.6}, style=thick, scale=0.5]
\node[label={[shift={(0.3,-0.8)}]$v_1$}] (V1) at (0,0){};
\node (V2) at (1.5,0){};
\node[label={[shift={(0.6,-1)}]$v_{h-1}$}] (V3) at (3,0){};
\node[label={[shift={(0.3,-0.8)}]$v_h$}] (V4) at (4.5,0){};
\draw (6,0) circle (1.5) node [fill = white, scale=1.75] {$C_{t}$}(V4);
\draw (-1.5,0) circle (1.5) node [fill = white, scale=1.75] {$C_{s}$}(V1);
\draw (V1) -- (V2);
\draw (V3) -- (V4);
\draw[dotted] (V2)-- (V3);
\end{tikzpicture}
\qquad \qquad \qquad
\begin{tikzpicture}[auto=center,every node/.style={circle, fill=black, scale=0.6}, style=thick, scale=0.5]
\node[label={[shift={(-0.3,-0.8)}]$u$}] (V1) at (0,0){};
\node (V2) at (1.5,0){};
\node (V3) at (3,0){};
\node (V4) at (4.5,0){};
\node (V5) at (6,0){};
\node[label={[shift={(0.3,-0.8)}]$v$}] (V6)  at (7.5,0){};
\node (V7) at (1.5,1.5){};
\node (V8) at (3,1.5){};
\node (V9) at (4.5,1.5){};
\node (V10) at (6,1.5){};
\node (V11) at (1.5,-1.5){};
\node (V12) at (3,-1.5){};
\node (V13) at (4.5,-1.5){};
\node (V14) at (6,-1.5){};
\draw (V1) -- (V2) -- (V3);
\draw (V4) -- (V5) -- (V6);
\draw (V1) -- (V7) -- (V8);
\draw (V9) -- (V10) -- (V6);
\draw (V1) -- (V11) -- (V12);
\draw (V13) -- (V14) -- (V6);
\draw[dotted] (V3)-- (V4) node[midway,above = 0.1pt,fill=white] {$P_{\ell+1}$};
\draw[dotted] (V8)-- (V9) node[midway,above = 0.1pt,fill=white] {$P_{k+1}$};
\draw[dotted] (V12)-- (V13) node[midway,above = 0.1pt,fill=white] {$P_{m+1}$};
\end{tikzpicture}

\caption{Left: a type-I bicyclic graph; right: a type-II bicyclic graph.}
\label{fig.bicyclic}
\end{figure}

A simple connected graph $G$ is said to be \defn{bicyclic} if the number of edges is equal to the number of vertices plus one. Let $\mathcal{B}_n$ denote the set of connected bicyclic graphs of order $n$. For any $G\in \mathcal{B}_n$, the unique minimal bicyclic subgraph $\tilde{G}$ falls into one of two cases: 
\begin{enumerate}
    \item $\tilde{G}$ is constructed from two disjoint cycles $C_s$ and $C_t$ by adding a path $P_h$ (with $h \geq 1$) from a vertex in $C_s$ to any vertex in $C_t$. In the case where $h=1$, we have two cycles that share a single vertex. 
    \item $\tilde{G}$ is constructed from three paths whose edge sets are disjoint, $P_{k+1}, P_{\ell+1}, P_{m+1}$, where $k,\ell,m \geq 1$ and at most one of them is equal to 1. The initial vertices of all three paths are identified as a single vertex $u$, and the terminal vertices are identified as a single vertex $v$. 
    In other words, $\tilde{G}$ can be regarded as consisting of two cycles $C_{\ell+k}$ and $C_{\ell+m}$ which share a common path of $P_{\ell+1}$. 
\end{enumerate}
For any $G\in \mathcal{B}_n$, we say that $G$ is a \defn{type-I} (resp. \defn{type-II}) bicyclic graph if its unique minimal bicyclic subgraph is of the form of case (1) (resp. case (2)) above. We illustrate these two types in \Cref{fig.bicyclic}. 

\begin{example} In \Cref{fig:typeI_ex} we illustrate a type-I bicyclic graph in $\mathcal{B}_{16}$. Its unique minimal bicyclic subgraph has cycles of size $4$ and $5$ connected by a path with $h = 3$.
\begin{figure}[hbt!]
\centering
\begin{tikzpicture}[auto=center,every node/.style={circle, fill=black, scale=0.45}, style=thick, scale=0.4] 

%%%%% Draw vertices %%%%%
\node (V1) at (0,0) {};
\node (V2) at (1.5,1.5) {};
\node (V3) at (3,0) {};
\node (V4) at (1.5,-1.5) {};
\node (V5) at (6,0) {};
\node (V6) at (7.5,1.5) {};
\node (V7) at (9,1.5) {};
\node (V8) at (9,-1.5) {};
\node (V9) at (7.5,-1.5) {};
\node (V10) at (4.5,0) {};
\node (V11) at (4,1.5) {};
\node (V12) at (5,1.5) {};
\node (V13) at (5,-1.5) {};
\node (V14) at (10.5,1.5) {};
\node (V15) at (10.5,0) {};
\node (V16) at (12,1.5) {};

%%%%% Draw edges %%%%%
\draw (V1)--(V2)--(V3)--(V4)--(V1);
\draw (V3)--(V10)--(V5);
\draw (V5)--(V6)--(V7)--(V8)--(V9)--(V5);
\draw (V12)--(V10)--(V11);
\draw (V5)--(V13);
\draw (V15)--(V14)--(V16);
\draw (V14)--(V6);

\end{tikzpicture}
\caption{A type-I bicyclic graph in $\mathcal{B}_{16}$}
\label{fig:typeI_ex}
\end{figure}
\end{example}

\begin{proposition}\label{prop.bicycles}
    Let $G$ be a bicyclic graph on $n$ vertices. 
    \begin{enumerate}
        \item If $G$ is a type-I bicyclic graph formed by cycles of length $s,t \geq 3$, then
        \[c_{(n)}=(s-1)(t-1).\]
        \item If $G$ is a type-II bicyclic graph formed by two cycles of length $s,t \geq 3$ which meet along a path $P_{\ell+1}$ for $\ell \geq 1$, then 
        \[c_{(n)}=(s-1)(t-1)-2\binom{\ell}{2}.\]   
        Note that when $\ell=1$, the formula reduces to $c_{(n)}=(s-1)(t-1)$.
    \end{enumerate}
\end{proposition}
\begin{proof}
First we consider case (1). We apply the DNC relation to an edge of the cycle containing $s$ vertices. The term corresponding to deletion of the edge contributes $t-1$ to the coefficient $c_{(n)}$, by \Cref{cor:cyc size from coeff}. Next, the dot-contraction term can be ignored since it produces a disconnected graph, for which $c_{(n)} = 0$ since all partitions from disconnected graphs have length at least $2$. We now proceed by induction on $s$ for the leaf-contraction term. When $s=3$, the leaf-contraction term also contributes $(t-1)$ yielding $2(t-1)$ total for $c_{(n)}$, as claimed. Now if $s>3$, the leaf-contraction term reduces us to the case of cycles of size $s-1$ and $t$, which contributes $(s-2)(t-1)$ by the inductive hypothesis. Finally, adding this to the contribution of deletion produces the result.

Next we consider case (2) where the cycles of size $s,t$ have $\ell \geq 1$ edges in common. We proceed with induction on $\ell$. For the base case, suppose $\ell=1$ and $e$ is the single shared edge. We apply the DNC relation to $e$. If $s=3$, deleting the edge $e$ yields a unicyclic graph with a cycle of size $2+(t-1)=t+1$, which by \Cref{cor:cyc size from coeff} contributes $t$ to $c_{(n)}$. If $t=3$, leaf-contracting $e$ would result in a tree which contributes $1=t-2$ to $c_{(n)}$ by \Cref{prop:treehook}. If $t \geq 4$, leaf-contracting $e$ would produce a unicyclic graph with a cycle of size $t-1$, which contributes $t-2$ to $c_{(n)}$.
We may again ignore the dot-contraction term. Thus by the DNC relation, $c_{(n)}=t+t-2=2(t-1)$. For $s>3$, applying DNC to $e$ now gives $s+t-3$ from deletion, and $(s-2)(t-2)$ from leaf-contraction by applying the formula for case (1). Adding these contributions yields the result. Supposing the result holds whenever the cycles have $\ell-1$ edges in common, we apply DNC to an edge of the common path. Deletion contributes $s+t-2\ell-1$, since we are left with one cycle of length $s+t-2\ell$. Leaf-contraction reduces to the $\ell-1$ case where our inductive hypothesis has contribution 
\[(s-2)(t-2)-2 \binom{\ell-1}{2}=st-2s-2t+4-(\ell-1)(\ell-2).\]
Adding these contributions, we have 
\[st-t-s+1-2(\ell-1)-(\ell-1)(\ell-2)=(s-1)(t-1)-\ell(\ell-1)\]
as desired.
\end{proof}

In what follows, we will show how to recover the sizes of the two cycles in a bicyclic graph. In order to do this we need another equation in addition to $(s-1)(t-1)=c_{(n)}$ in type-I, and we need two additional equations in type-II as we have three unknowns. We begin with a proposition. 

\begin{proposition}{\cite[Lemma 4.19]{OrellanaTom25}}\label{prop:unicyclic_coeff_sum}
Let $j\ge 1$, $G$ be a unicyclic graph with a cycle of length $c \ge 3$ and chromatic symmetric function $\mathbf{X}_G = \sum_{\lambda\vdash n} c_\lambda \mf{st}_\lambda$. Then, \[
\sum_{\substack {\lambda\vdash n \\ \ell(\lambda) = j}} c_\lambda = (-1)^{j-1}\binom{c-1}{j}. 
\]
In particular, these sums are zero whenever $j\geq c$.
\end{proposition}
We have some similar results for bicyclic graphs. 

\begin{proposition}\label{prop:bicyclicI_coeff_sum}
Let $G$ be a connected bicyclic graph on $n$ vertices of type-I containing cycles of length $t\ge s \ge 3$. 
Let $\textbf{X}_{G} = \sum_{\lambda \vdash n} c_{\lambda} \mathfrak{st}_{\lambda}$. Then for $j \ge 1$,
$$
    \sum_{\lambda \vdash n, \ \ell(\lambda)=j} c_{\lambda} 
    = 
    (-1)^{j-1} \sum_{i=1}^{j} \tbinom{s-1}{i} \tbinom{t-1}{j-i+1}
    =(-1)^{j-1}\left[\tbinom{s+t-2}{j+1}-\tbinom{t-1}{j+1}-\tbinom{s-1}{j+1}\right]
    \,.
$$   
In particular, the sum vanishes when $j\geq s+t-2$ and is non-zero for $1\leq j \leq s+t-3$.  
\end{proposition}

\begin{proof}
We induct on the length of the cycle $s$. 
For the base case, suppose $s=3$. We apply the DNC relation to an edge $e$ on the $s$-cycle, we have that $G\setminus e$ and $G \odot e$ are both connected unicyclic graphs on $n$ vertices containing a $t$-cycle, while $(G \odot e) \setminus \ell_e$ is the disjoint union of a connected unicyclic graph on $n-1$ vertices and an isolated vertex. 
Hence, by \cref{prop:unicyclic_coeff_sum}, $\textbf{X}_{G \setminus e}$ and $\textbf{X}_{G \odot e}$ each contribute $(-1)^{j-1} \binom{t-1}{j}$ to the sum, and $\textbf{X}_{(G \odot e)\setminus \ell_e}$ contributes $0$ to the sum when $j=1$, and contributes $(-1)^{j-2}\binom{t-1}{j-1}$ when $j >1$.
    Therefore, by the DNC relation, $\sum_{\lambda \vdash n,\  \ell(\lambda)=j} c_{\lambda} 
        = (-1)^{j-1} 2\tbinom{t-1}{j} = (-1)^{j-1}\tbinom{3-1}{1} \tbinom{t-1}{j}$ when $j=1$ and when $j>1$,
    $$
    \begin{aligned}
        \sum_{\lambda \vdash n,\  \ell(\lambda)=j} c_{\lambda} 
        = 
        (-1)^{j-1}\tbinom{t-1}{j}
        -(-1)^{j-2} \tbinom{t-1}{j-1}
        +(-1)^{j-1} \tbinom{t-1}{j}
        &=(-1)^{j-1} \left[ 2 \tbinom{t-1}{j} + \tbinom{t-1}{j-1}\right]
        \\
        &
        =
        (-1)^{j-1} \left[ \tbinom{3-1}{1}\tbinom{t-1}{j} + \tbinom{3-1}{2} \tbinom{t-1}{j-1}\right]
        \,.
    \end{aligned}
    $$
This completes the base case. 
For the inductive step, assume that the statement holds for all such bicyclic graphs of type-I with cycles of size $t,z-1 \ge 3$. 
Let $G$ be a bicyclic graph on $n$ vertices with cycles of size $t,z \ge 3$ sharing at most one vertex. 
By applying the DNC relation to an edge $e$ on the $z$-cycle, we have that $G \setminus e$ is a connected unicyclic graph on $n$ vertices containing a $t$-cycle; $(G \odot e)\setminus  \ell_e$ is the disjoint union of a bicyclic graph (of type-I) on $n-1$ vertices containing cycles of sizes $t$ and $z-1$ and an isolated vertex; and $G \odot e$ is a bicyclic graph on $n$ vertices (of type-I) containing cycles of size $t$ and $z-1$. 
Hence, by \cref{prop:unicyclic_coeff_sum}, $\textbf{X}_{G \setminus e}$ contributes $(-1)^{j-1}\binom{t-1}{j}$ to the sum;
$\textbf{X}_{(G \odot e)\setminus \ell_e}$ contributes  
$(-1)^{j-2} \sum_{i=1}^{j-1} \binom{z-2}{i} \binom{t-1}{j-i}$ by the inductive hypothesis;
and, again by the inductive hypothesis, 
$\textbf{X}_{G \odot e}$ contributes $(-1)^{j-1} \sum_{i=1}^{j} \binom{z-2}{i} \binom{t-1}{j-i+1}$. 
Therefore, we have that
  $$
    \begin{aligned}
        \sum_{\lambda \vdash n,\  \ell(\lambda)=j} c_{\lambda} 
        &= 
        (-1)^{j-1} \tbinom{t-1}{j} - (-1)^{j-2} \sum_{i=1}^{j-1} \tbinom{z-2}{i}\tbinom{t-1}{j-i} + (-1)^{j-1} \sum_{i=1}^{j} \tbinom{z-2}{i} \tbinom{t-1}{j-i+1}\\
        &=
        (-1)^{j-1} \left[\tbinom{t-1}{j} + \sum_{i=1}^{j-1} \tbinom{z-2}{i}\tbinom{t-1}{j-i} + \sum_{i=1}^{j} \tbinom{z-2}{i} \tbinom{t-1}{j-i+1}\right]
        \\
        &=
        (-1)^{j-1} \left[\tbinom{t-1}{j} 
        + 
        \sum_{i=1}^{j-2} \tbinom{z-2}{i}\tbinom{t-1}{j-i} 
        + 
        \tbinom{z-2}{j-1}\tbinom{t-1}{1}
        +
        \tbinom{z-2}{1}\tbinom{t-1}{j}
        + 
        \sum_{i=2}^{j-1} \tbinom{z-2}{i} \tbinom{t-1}{j-i+1} + \tbinom{z-2}{j}\tbinom{t-1}{1}\right]
        \,.
    \end{aligned}
    $$
Reindex the last sum, combine with terms of the second, and apply Pascal's rule to get
    $$
    \begin{aligned}
        \sum_{\lambda \vdash n,\  \ell(\lambda)=j} c_{\lambda} 
        &= 
        (-1)^{j-1} \left[\tbinom{z-1}{1}\tbinom{t-1}{j} 
        + \sum_{i=1}^{j-2} 
        \left(\tbinom{z-2}{i}
        +\tbinom{z-2}{i+1} \right)\tbinom{t-1}{j-i} 
        + \tbinom{z-2}{j-1}\tbinom{t-1}{1}
        + \tbinom{z-2}{j}\tbinom{t-1}{1}
        \right]
        \\
        &=
        (-1)^{j-1} \left[\tbinom{z-1}{1}\tbinom{t-1}{j} 
        + \sum_{i=1}^{j-2} 
        \tbinom{z-1}{i+1} \tbinom{t-1}{j-i} 
        + \tbinom{z-1}{j}\tbinom{t-1}{1}
        \right]
        \\
        &=
        (-1)^{j-1} \left[\sum_{i=1}^{j} \tbinom{z-1}{i}\tbinom{t-1}{j-i+1}
        \right]
        \,.
    \end{aligned}
    $$
The result follows from applying Vandermonde's identity. 

It is easy to see that if $j\geq s+t-2$, then all binomial coefficients ${s+t-2\choose j+1}$, ${s-1\choose j+1}$, and ${t-1\choose j+1}$ are zero. Thus, the sums vanish for these values of $j$. Now, we show that the summation is strictly positive for $1\leq j \leq s+t-3$. Since
\[\tbinom{s+t-2}{j+1} = \sum_{i=0}^{j+1} \tbinom{s-1}{i} \tbinom{t-1}{j-i+1} = \sum_{i=1}^{j} \tbinom{s-1}{i} \tbinom{t-1}{j-i+1}+\tbinom{t-1}{j+1}+\tbinom{s-1}{j+1},\] then $\tbinom{s+t-2}{j+1}>\tbinom{t-1}{j+1}+\tbinom{s-1}{j+1}$ if there is at least one $i$ such that $ \tbinom{s-1}{i} \tbinom{t-1}{j-i+1}>0$. 
This happens when there is an $i$ in the summation such that $1\leq i\leq s-1$ and $1\leq j-i+1\leq t-1$, which is equivalent to $1\leq i \leq s-1$ and $j-t+2\leq i \leq j$. Thus, there is an $i$ in both intervals if $1\leq j$ and $j-t+2\leq s-1$, in other words $1\leq j \leq s+t-3$.
\end{proof}

\begin{proposition} \label{prop:bicyclicIIcoeff}
Let $G$ be a connected bicyclic graph of type-II, formed by two cycles of length $s,t \geq 3$ which meet along a path $P_{\ell+1}$ for $\ell \geq 1$, and with chromatic symmetric function $\mathbf{X}_G = \sum_{\lambda\vdash n} c_\lambda \mf{st}_\lambda$. Then
\[
\sum_{\substack {\lambda\vdash n \\ \ell(\lambda) = j}}c_\lambda = (-1)^{j-1}\left[\sum_{i = 1}^\ell \left[\binom{s+t-2(\ell+1)+i}{j}\right] + \binom{s+t-\ell-2}{j+1} - \binom{s-1}{j+1} - \binom{t-1}{j+1}\right].
\]

In particular, the sum vanishes when $j \ge s+t-\ell-1$ and is non-zero for $1\leq j< s+t-\ell -1$. We also note that the formula reduces to that of a type-I bicyclic graph when $\ell = 1$, and, although type-II bicyclic graphs assume $\ell \ge 1$, setting $\ell = 0$ also agrees with the type-I formula.
\end{proposition}

\begin{proof}
We proceed by induction on the length of the common path $P_{\ell+1}$ and assume without loss of generality that $s \le t$. For the base case, suppose $\ell = 1$. If $s = 3$, we let $e$ be an edge of the cycle containing $s$ vertices that is not the single edge in $P_{\ell+1}$. We apply the DNC relation to $e$ to obtain $\mathbf{X}_G =\mathbf{X}_{G\setminus e} -\mathbf{X}_{(G\odot e)\setminus \ell_e} + \mathbf{X}_{G\odot e}$ and consider $\sum_{\substack {\lambda\vdash n \\ \ell(\lambda) = j}} c_\lambda$ for each of the three graphs coming from the DNC relation. Notice that the deletion and leaf-contraction of $e$ each result in a connected unicyclic graph with a cycle of length $t$ and the dot-contraction of $e$ results in a unicyclic graph with a cycle of length $t$ together with an isolated vertex. For $j = 1$, $\mathbf{X}_{G\setminus e}$ and $\mathbf{X}_{G\odot e}$ each contribute $t-1$ to the sum by \Cref{cor:cyc size from coeff} while $\mathbf{X}_{(G\odot e)\setminus \ell_e}$ contributes $0$ to the sum since it is disconnected. The DNC relation and Pascal's rule show
\begin{align*}
    c_{(n)} &= 2t - 2 = \tbinom{t}{1} +\tbinom{t -1}{1} - \tbinom{2}{2}\\ &= \tbinom{t}{1} + \tbinom{t}{2} - \tbinom{2}{2} - \tbinom{t-1}{2}\\ &= (-1)^0 \tbinom{3+t-2(2)+1}{1} + \tbinom{3+t-1-2}{2} - \tbinom{3-1}{2} - \tbinom{t-1}{2}. 
\end{align*}
For $j \ge 2$, \Cref{prop:unicyclic_coeff_sum} states that $\mathbf{X}_{G\setminus e}$ and $\mathbf{X}_{G\odot e}$ each contribute $(-1)^{j-1}\binom{t-1}{j}$ to the sum. With $\mathbf{X}_{(G\odot e)\setminus \ell_e}$, all partitions of length $j$ result from multiplying those of length $j-1$ by $\mf{st}_{(1)}$, and therefore by \Cref{prop:unicyclic_coeff_sum}, $\mathbf{X}_{(G\odot e)\setminus \ell_e}$ contributes $(-1)^{j-2}\binom{t-1}{j-1}$ to the sum.
Since $j \ge 2$, it follows that $\tbinom{2}{j+1}=0$. The DNC relation and Pascal's rule show
\begin{align*}
    \sum_{\substack {\lambda\vdash n \\ \ell(\lambda) = j}} c_\lambda 
    &= (-1)^{j-1} \left[\tbinom{t-1}{j} +\tbinom{t -1}{j-1} +\tbinom{t-1}{j}- \tbinom{2}{j+1}\right]\\ 
    &= (-1)^{j-1}\left[\tbinom{t}{j} + \tbinom{t-1}{j} - \tbinom{2}{j+1}\right]\\ 
    &= (-1)^{j-1}\left[\tbinom{t}{j} + \tbinom{t}{j+1} - \tbinom{2}{j+1} - \tbinom{t-1}{j+1}\right]\\ 
    &= (-1)^{j-1} \left[\tbinom{3+t-2(2)+1}{j} + \tbinom{3+t-1-2}{j+1} - \tbinom{3-1}{j+1} - \tbinom{t-1}{j+1}\right]. 
\end{align*}
If $s \ge 4$, we let $e$ be the single shared edge. Notice that $G\setminus e$ is a connected unicyclic graph with cycle length $s+t - 2$, $(G\odot e)\setminus \ell_e$ is a type-I graph with cycles of length $s-1$ and $t-1$ together with an isolated vertex, and $(G\odot e)$ is a connected type-I graph with cycles of size $s-1$ and $t-1$. By \Cref{prop:unicyclic_coeff_sum}, $\mathbf{X}_{G\setminus e}$ contributes $(-1)^{j-1}\binom{s+t-3}{j}$ to the sum. By \Cref{prop:bicyclicI_coeff_sum}, $\mathbf{X}_{(G\odot e)\setminus \ell_e}$ contributes $(-1)^{j-2} \left[\binom{s+t-4}{j}-\binom{s-2}{j} - \binom{t-2}{j}\right]$ and $\mathbf{X}_{G\odot e}$ contributes $(-1)^{j-1} \left[\binom{s+t-4}{j+1}-\binom{s-2}{j+1} - \binom{t-2}{j+1}\right]$. Using these contributions in the DNC relation and Pascal's rule, we have 
\begin{align*}
    \sum_{\substack {\lambda\vdash n \\ \ell(\lambda) = j}} c_\lambda 
    &= (-1)^{j-1}\left[\tbinom{s+t-3}{j} + \tbinom{s+t-4}{j}-\tbinom{s-2}{j} - \tbinom{t-2}{j} + \tbinom{s+t-4}{j+1}-\tbinom{s-2}{j+1} - \tbinom{t-2}{j+1}\right]\\
    &= (-1)^{j-1}\left[\tbinom{s+t-3}{j} + \tbinom{s+t-3}{j+1}-\tbinom{s-1}{j+1} - \tbinom{t-1}{j+1}\right].
\end{align*}

For the inductive step, suppose the statement holds for $\ell = w \ge 1$ and let $G$ be as in the statement with $\ell = w+1$. Observe that we may assume that $s, t \ge 4$ since if there were a $3$-cycle, we could rearrange the paths so that the common path has length $1$, which is covered by the base case. Let $e$ be an edge of the common path, $P_{w+2}$. As above, we apply the DNC relation on $e$ and analyze the contributions to the sum from each of the three graphs. The deletion of $e$ results in a unicyclic graph with cycle length $s+t - 2(w+1)$ and thus $\mathbf{X}_{G\setminus e}$ contributes $(-1)^{j-1}\binom{s+t-2(w+1)-1}{j}$ to the sum. The dot-contraction of $e$ results in a type-II bicyclic graph with cycles of size $s-1$ and $t-1$ that meets along a path of length $w$ together with an isolated vertex and thus by assumption contributes $(-1)^{j-2}\left(\sum_{i = 1}^w \left[\binom{(s-1)+(t-1)-2(w+1)+i}{j-1}\right] + \binom{(s-1)+(t-1)-w-2}{j} - \binom{s-2}{j} - \binom{t-2}{j}\right)$. Similarly, the leaf-contraction of $e$ contributes $(-1)^{j-1}\left(\sum_{i = 1}^w \left[\binom{(s-1)+(t-1)-2(w+1)+i}{j}\right] + \binom{(s-1)+(t-1)-w-2}{j+1} - \binom{s-2}{j+1} - \binom{t-2}{j+1}\right)$. Thus by the DNC relation and Pascal's rule, 
\allowdisplaybreaks
\begin{align*}
    \sum_{\substack {\lambda\vdash n \\ \ell(\lambda) = j}} c_\lambda 
    &= (-1)^{j-1}\bigg[\tbinom{s+t-2(w+1)-1}{j}\\
    & \qquad \qquad \qquad + \sum_{i = 1}^w \left[\tbinom{s +t-2(w+2)+i}{j-1}\right] + \tbinom{s+t-(w+1)-3}{j} - \tbinom{s-2}{j} - \tbinom{t-2}{j} \\
    & \qquad \qquad \qquad + \sum_{i = 1}^w \left[\tbinom{s+t-2(w+2)+i}{j}\right] + \tbinom{s+t-(w+1)-3}{j+1} - \tbinom{s-2}{j+1} - \tbinom{t-2}{j+1}\bigg]\\
    &=(-1)^{j-1}\left[\tbinom{s+t-2(w+2)+1}{j} + \sum_{i = 2}^{w+1} \left[\tbinom{s +t-2(w+2)+i}{j}\right] + \tbinom{s+t-(w+1)-2}{j+1} - \tbinom{s-1}{j+1} - \tbinom{t-1}{j+1}\right],
\end{align*}
which after combining the first two terms as a single sum, we get the desired result. 

It is easy to see that all binomials in the statement of the proposition vanish when $j\geq s+t-\ell -1$.
If $\ell =1$, then the formula reduces to type-I, hence the non-vanishing is proven in Proposition \ref{prop:bicyclicI_coeff_sum}. In fact, if any of $k,m,$ or $\ell$ are equal to 1, the formula reduces to that of type-I. Hence, we assume that $k,m, \ell \geq 2$ (where $k, m, \ell$ are as in our definition of bicyclic graphs, see Figure \ref{fig.bicyclic}).  Next we show that when $1\leq j \leq s+t-\ell-2$, $\sum_{\substack {\lambda\vdash n \\ \ell(\lambda) = j}} c_\lambda \neq 0$. 
We claim that 
\begin{equation}\label{eq:WTS}
\sum_{i = 1}^\ell \binom{s+t-2(\ell+1)+i}{j} + \binom{s+t-\ell-2}{j+1} - \binom{s-1}{j+1} - \binom{t-1}{j+1}
\end{equation}
is greater than zero if $1\leq j \leq s+t-\ell - 2$.
Notice that it is greater than zero whenever $j\geq t-1$ since we assume that $s\leq t$, and for $j\geq t-1$ the last two subtracted binomials are zero.  It remains to prove that expression (\ref{eq:WTS}) is positive for $1\leq j \leq t-2$.
To show this, we use the following identity that is a direct consequence of Pascal's relation for any two positive integers $a>b$ we have 
\begin{equation} \label{eq:telescope}
\binom{a}{j+1} - \binom{b}{j+1} = \sum_{i = b}^{a-1} \binom{i}{j}.
\end{equation}
We also use the fact that $\binom{a}{i} > \binom{b}{i}$ whenever $a> b\geq 0$ and $1\leq i \leq a$. 
By Equation (\ref{eq:telescope}) with $a=s+t-\ell -2$ and $b = t-1$, we get
\[
\binom{s+t-\ell -2}{j+1}-\binom{t-1}{j+1} = \sum_{i=t-1}^{s+t-\ell -3}\binom{i}{j} \geq 0
\]
Notice that this means that since $\binom{s-1}{j+1} =0$ for $j\geq s-1$, we actually only need to show that expression (\ref{eq:WTS}) is positive for $1\leq j \leq s-2$.
We use \cref{eq:telescope} with $a=s-1$ and $b=j$, to show that 
\begin{equation}
\binom{s-1}{j+1} = \sum_{i = j}^{s-2}\binom{i}{j}.
\end{equation}
Substituting into expression (\ref{eq:WTS}) we get 
\begin{align}
    \abs{\sum_{\substack {\lambda\vdash n \\ \ell(\lambda) = j}} c_\lambda} &= \sum_{i = 1}^\ell \binom{s+t-2(\ell+1)+i}{j} +\sum_{i=t-1}^{s+t-\ell-3}\binom{i}{j} - \sum_{i=j}^{s-2}\binom{i}{j} \nonumber \\
&=\sum_{i = 1}^\ell \binom{s+t-2(\ell+1)+i}{j} +\left(\sum_{i=t-1}^{s+t-\ell-3}\binom{i}{j} - \sum_{i=l}^{s-2}\binom{i}{j}\right)-\sum_{i=j}^{\ell-1}\binom{i}{j}
.\label{eq:step1}
\end{align}
Examining the middle expression of \eqref{eq:step1} and reindexing, we have
\begin{equation}\label{eq:step3} \sum_{i=t-1}^{s+t-\ell-3}\binom{i}{j} - \sum_{i=\ell}^{s-2}\binom{i}{j} =\sum_{i=0}^{s-\ell-2}\left[\binom{t-1+i}{j}-\binom{\ell+i}{j}\right]> 0\,,
\end{equation}
as $t-1 = (\ell + m) - 1 > \ell$ and $s - \ell -2 = (\ell + k) - \ell + 2 \geq 0$.

This implies that it suffices to show
\begin{equation}\label{eq:step2}\sum_{i = 1}^\ell \binom{s+t-2(\ell+1)+i}{j}-\sum_{i=j}^{\ell-1}\binom{i}{j} \geq 0.\end{equation}
If $\ell\leq j$ then the subtracted summation of inequality (\ref{eq:step2}) is empty and we are done. Otherwise, if $\ell>j$, we can reindex both summations and get 
\begin{align*}
    \sum_{i = 1}^\ell \binom{s+t-\ell-1-i}{j}-\sum_{i=1}^{\ell-j}\binom{\ell-i}{j}
    &\geq \sum_{i = 1}^\ell \binom{t-1-i}{j}-\sum_{i=1}^{\ell-j}\binom{\ell-i}{j}\\
    &= \sum_{i=\ell-j+1}^\ell \binom{t-1-i}{j}+\left(\sum_{i = 1}^{\ell-j} \left[\binom{t-1-i}{j}-\binom{\ell-i}{j}\right]\right) \geq 0,
\end{align*}
where the first inequality follows from $\ell \leq s-1$, and last inequality follows from $\ell\leq t-1$. So, the summation in parenthesis is greater or equal to zero. By inequality (\ref{eq:step3}) this finishes the proof.
\end{proof}

Recall that the \emph{girth} of a graph $G$ is defined as the length of the smallest cycle. To show that we can recover the cycle sizes of bicyclic graphs, we will make use of the following result.

\begin{proposition}[{\cite[Prop. 3]{martin2007distinguishingtreeschromaticsymmetric}}]\label{prop:MMWgirth}
    The girth $g$ of a graph $G$ is recoverable from $\mathbf{X}_G$. 
\end{proposition}

Notice that in type-II, the formulas of Propositions~\ref{prop.bicycles}(2) and \ref{prop:bicyclicIIcoeff} are independent of the choice of $s$ and $t$. On the other hand, we can always assume that $\ell \leq k  \leq m$ and that $s = \ell + k$ is the girth of the graph so that $t=\ell+m$ (see Figure \ref{fig.bicyclic}). That is, the path common to the two cycles has length $\ell \leq s/2<s$. We will make this assumption to facilitate the proof of our main result for this section.

\begin{theorem}
\label{thm.bicyclic.recover}
For any bicyclic graph $G$ with cycle sizes $s\leq t$ meeting along $\ell<s/2$ edges, we can recover $s$ and $t$ from $\mathbf{X}_G$. Furthermore, if $\ell\geq 2$ we can detect that $G$ is of type-II and recover $\ell$.
\end{theorem}

\begin{proof}
    First note that the girth of a bicyclic graph is equal to the size of the smaller cycle, $s=\ell+k$, so this is immediately recovered from \cref{prop:MMWgirth}. Next we show how we can distinguish between the cases $\ell\geq 2$ and $\ell \leq 1$. Let
    \[\displaystyle{j^*:=\min \{j\in [n] \mid \sum_{\substack{\la\vdash n \\ \ell(\la)=j}}c_\la=0\} },\]
    and 
    \[\hat{t}:=\frac{c_{(n)}}{s-1}+1.\]
    Observe that if $G$ is of type-I, i.e., $\ell =0$, or is of type-II with $\ell =1$, then $t=\hat{t}$,  while if $G$ is of type-II with $\ell \geq 2$ then the size $t=\ell+m$ of the other cycle satisfies
    \[t=\frac{c_{(n)}+2\binom{\ell}{2}}{s-1}+1=\hat{t}+\frac{\ell(\ell-1)}{s-1}>\hat{t}.\]
Observe that if $\hat{t}$ is not integral, then we may immediately conclude $G$ is of type-II and that $\ell\geq 2$. Otherwise, it must be the case that $s-1 \mid \ell(\ell-1)$, and we proceed as follows.

By \cref{prop:bicyclicI_coeff_sum} and \cref{prop:bicyclicIIcoeff}, we have that if $G$ is type-I or $\ell =1$ then $j^*=s+\hat{t}-2$ while if $G$ is of type-II then $j^*=s+t-\ell-1$. We claim that $\ell \geq 2$ exactly when $s+\hat{t}-2\neq j^*$.
This will follow from the inequality 
\[\hat{t}-1\geq t-\ell\]
which we claim is strict for $\ell\geq 2$. Observe that if $\ell\geq 2$ we have
\begin{align*}
    t-\ell&=\frac{c_{(n)}+\ell(\ell-1)}{s-1}+1-\ell\\
    &=\hat{t}-1+\left(1-\frac{\ell(s-\ell)}{s-1}\right).
\end{align*} 
The claim now follows from the basic inequality $\ell(s-\ell)> s-1$ for $1<\ell<s-1$, with equality at the endpoints. Note that with our assumptions that $\ell \leq k$ and $s = \ell + k$ (see the paragraph preceding the theorem statement), $\ell$ is always at most $s/2$ and when $\ell \geq 2$, we have $s>3$.

Having shown that we can distinguish the case $\ell\leq 1$ (including type-I where $\ell=0$) from type-II with $\ell \geq 2$, we can conclude that if $G$ is type-I or $\ell=1$, we have $t=\hat{t}=j^*+2-s$. On the other hand, if $G$ is type-II we can obtain $\ell$ by identifying the equations $t=\frac{c_{(n)}+2\binom{\ell}{2}}{s-1}+1=j^*-s+\ell+1$ and solving for $\ell$ using the quadratic equation. This produces solutions
\[ \ell =\frac{1}{2}(s \pm \sqrt{s^2-4((s-j^*)(s-1)+c_{(n)})}).\]
Then substituting for $j^*$ and $c_{(n)}$, we see the discriminant gives
\begin{align*}
    s^2-4((s-1)(\ell+1-t)+(s-1)(t-1)-\ell(\ell-1)) &=s^2-4((s-1)(\ell)-\ell(\ell-1)) \\ 
    &=s^2-4\ell(s-\ell)\\
    &=(s-2\ell)^2.
\end{align*}
This shows that the solutions to the quadratic are exactly the natural number values $\ell$ and $k=s-\ell$, as desired. The smaller of these solutions in our convention is $\ell$, which can then be used with $j^*$ to recover $t=j^*-s+\ell+1$. 
\end{proof}

We remark that it is impossible to distinguish between type-I bicyclic graphs and the $\ell=1$ case of type-II bicyclic graphs. Indeed, the first example provided by Stanley of non-isomorphic graphs with the same chromatic symmetric function exhibits this difficulty \cite{stanley_1995} (see~\Cref{fig: type_i_equals_type_ii}).

\begin{example}
\begin{figure}[hbt!]
    \centering
    \begin{tikzpicture}[auto=center,every node/.style={circle, fill=black, scale=0.45}, style=thick, scale=0.5] 
        %%%%% TREE 1 %%%%%
    %%%%% Draw vertices %%%%%
    \filldraw[black] (0, 0) coordinate (A0) circle (4pt) node{};
    \filldraw[black] (-1, 1) coordinate (A1) circle (4pt) node{};
    \filldraw[black] (-1, -1) coordinate (A2) circle (4pt) node{};
    \filldraw[black] (1, 1) coordinate (A3) circle (4pt) node{};
    \filldraw[black] (1, -1) coordinate (A4) circle (4pt) node{};
    %%%%% Draw edges %%%%%
    \draw(A1) -- (A0) -- (A4) -- (A3) --  (A0) -- (A2) -- (A1);
        %%%%% TREE 2 %%%%%
    %%%%% Draw vertices %%%%%
    \filldraw[black] (4, 1) coordinate (B0) circle (4pt) node{};
    \filldraw[black] (4, -1) coordinate (B1) circle (4pt) node{};
    \filldraw[black] (6, -1) coordinate (B2) circle (4pt) node{};
    \filldraw[black] (6, 1) coordinate (B3) circle (4pt) node{};
    \filldraw[black] (8, 0) coordinate (B4) circle (4pt) node{};
    %%%%% Draw edges %%%%%
    \draw(B0) -- (B1) -- (B2) -- (B3) --  (B0) -- (B2);
    \draw(B3) -- (B4);
    \end{tikzpicture}
     \caption{Two bicyclic graphs with identical chromatic symmetric functions}
     \label{fig: type_i_equals_type_ii}
\end{figure}
Consider the two graphs in~\Cref{fig: type_i_equals_type_ii}. They both have chromatic symmetric function:
\[
\mathbf{X}_G=4\mf{st}_{(5)}-8\mf{st}_{(4,1)}+4\mf{st}_{(3,2)}+4\mf{st}_{(3,1,1)}-3\mf{st}_{(2,2,1)}.
\]
Now suppose we only know that $\mathbf{X}_G$ comes from a certain bicyclic graph $G$ with cycle lengths $s,t$ with $s\leq t$. We can recover the number of vertices $n = 5$ from the degree. Next we use \cite[Prop. 3]{martin2007distinguishingtreeschromaticsymmetric} to find $s = 3$.
Then we may follow the instructions from the proof of \Cref{thm.bicyclic.recover}:
\begin{enumerate}
\item $j^*:=\min \{j\in [n] \mid \sum_{\substack{\la\vdash n \\ \ell(\la)=j}}c_\la=0\} = 4$\,,
\item $\hat{t}:=\frac{c_{(n)}}{s-1}+1 = \frac{4}{3-1}+1=3$\,,
\item $s+\hat{t}-2 = 3+3-2 = 4 = j^*$\,.
\end{enumerate}
Therefore we must have $\ell < 2$ and $t = \hat{t} = 3$.
\end{example}

On the other hand, we have the following corollary.
\begin{corollary}
    No bicyclic graph of type-I or of type-II with $\min\{k,\ell, m\}=1$ can have the same CSF as a type-II bicyclic graph with $\min\{k,\ell, m\}\geq 2$ where $k, \ell, m$ are as in our definition of bicyclic graphs.
\end{corollary}

\section{Examples and data}
\label{sec: data}
It is known that a graph is bipartite if and only if the graph contains no cycles of odd length. This implies that every tree is a bipartite graph. In \cite{ORELLANA20141} the authors construct an infinite number of unicyclic graphs with a 3-cycle with the same CSF. A natural question to then ask is whether the bipartite property of trees plays a role in the tree isomorphism problem. That is, one might ask: If $G_1$ and $G_2$ are non-isomorphic bipartite graphs, is it always true that $\mathbf{X}_{G_1}\neq\mathbf{X}_{G_2}$?
We exhibit in \Cref{fig: n=12,fig:n=13} examples showing that this is not the case. We note that these are the only pairs of non-isomorphic connected 4-unicyclic graphs up to $n=17$ vertices.

 \begin{figure}[hbt!]
    \centering
    \begin{tikzpicture}[auto=center,every node/.style={circle, fill=black, scale=0.45}, style=thick, scale=0.4] 
        %%%%% TREE 1 %%%%%
    %%%%% Draw vertices %%%%%
    \filldraw[black] (0, 1) coordinate (A0) circle (4pt) node{};
    \filldraw[black] (0, 0) coordinate (A1) circle (4pt) node{};
    \filldraw[black] (1, 0) coordinate (A2) circle (4pt) node{};
    \filldraw[black] (1, 1) coordinate (A3) circle (4pt) node{};
    \filldraw[black] (2, -1) coordinate (A4) circle (4pt) node{};
    \filldraw[black] (2, 2) coordinate (A5) circle (4pt) node{};
    \filldraw[black] (2, 0) coordinate (A6) circle (4pt) node{};
    \filldraw[black] (-1, -2) coordinate (A7) circle (4pt) node{};    \filldraw[black] (-2, -3) coordinate (A8) circle (4pt) node{};
    \filldraw[black] (0, -1) coordinate (A9) circle (4pt) node{};    \filldraw[black] (0, -2) coordinate (A10) circle (4pt) node{};
     \filldraw[black] (3, -2) coordinate (A12) circle (4pt) node{};
    
    %%%%% Draw edges %%%%%
    \draw(A1) -- (A2) -- (A3) -- (A0) -- (A1);
    \draw(A3) -- (A5);
    \draw(A1) -- (A9) -- (A10);
    \draw(A9) -- (A7) -- (A8);
    \draw(A2) -- (A6);
    \draw(A2)-- (A4) -- (A12);

        %%%%% TREE 2 %%%%%
    %%%%% Draw vertices %%%%%
    \filldraw[black, ] (10, 1) coordinate (A0) circle (4pt) node{};
    \filldraw[black] (10, 0) coordinate (A1) circle (4pt) node{};
    \filldraw[black] (11, 0) coordinate (A2) circle (4pt) node{};
    \filldraw[black] (11, 1) coordinate (A3) circle (4pt) node{};
    \filldraw[black] (12, 2) coordinate (A5) circle (4pt) node{};
    \filldraw[black] (12, 0) coordinate (A6) circle (4pt) node{};
    \filldraw[black] (9, -2) coordinate (A7) circle (4pt) node{};    \filldraw[black] (8, -3) coordinate (A8) circle (4pt) node{};
    \filldraw[black] (10, -1) coordinate (A9) circle (4pt) node{};    \filldraw[black] (10, -2) coordinate (A10) circle (4pt) node{};
    \filldraw[black] (9, 0) coordinate (A11) circle (4pt) node{};    \filldraw[black] (13, 3) coordinate (A12) circle (4pt) node{};
    
    %%%%% Draw edges %%%%%
    \draw(A1) -- (A2)--(A3)--(A0)--(A1) --(A11);
    \draw(A3)--(A5)--(A12);
    \draw(A1)--(A9)--(A10);
    \draw(A9)--(A7)--(A8);
    \draw(A2)--(A6);
    
    \end{tikzpicture}
     \caption{Two bipartite graphs with the same CSF on $12$ vertices.}
     \label{fig: n=12}
\end{figure}
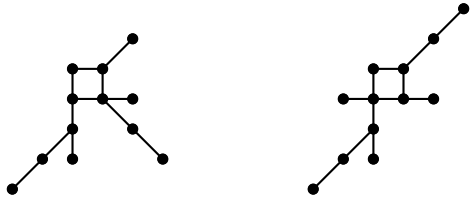

 %%% 13 vertex graphs
\begin{figure}[hbt!]
\centering    
\begin{tikzpicture}[auto=center,every node/.style={circle, fill=black, scale=0.45}, style=thick, scale=0.4] 
    %%%%% TREE 1 %%%%%
%%%%% Draw vertices %%%%%
\filldraw[black] (0, 1) coordinate (A0) circle (4pt) node{};
\filldraw[black] (0, 0) coordinate (A1) circle (4pt) node{};
\filldraw[black] (1, 0) coordinate (A2) circle (4pt) node{};
\filldraw[black] (1, 1) coordinate (A3) circle (4pt) node{};
\filldraw[black] (2, -1) coordinate (A4) circle (4pt) node{};
\filldraw[black] (2, 2) coordinate (A5) circle (4pt) node{};
\filldraw[black] (-1, 0) coordinate (A6) circle (4pt) node{};
\filldraw[black] (3, -2) coordinate (A7) circle (4pt) node{};
\filldraw[black] (4, -3) coordinate (A8) circle (4pt) node{};    
\filldraw[black] (5, -4) coordinate (A9) circle (4pt) node{};
\filldraw[black] (0, -1) coordinate (A10) circle (4pt) node{};    
\filldraw[black] (0, -2) coordinate (A11) circle (4pt) node{};
 \filldraw[black] (4, -2) coordinate (A12) circle (4pt) node{};

%%%%% Draw edges %%%%%
\draw(A1) -- (A2) -- (A3) -- (A0) -- (A1);
\draw(A3) -- (A5);
\draw(A1) -- (A10) -- (A11);
\draw (A7) -- (A8) -- (A9);
\draw(A1) -- (A6);
\draw(A2)-- (A4) -- (A7) -- (A12);
    %%%%% TREE 2 %%%%%
        %%%%% Draw vertices %%%%%
\filldraw[black] (10, 1) coordinate (A0) circle (4pt) node{};
\filldraw[black] (10, 0) coordinate (A1) circle (4pt) node{};
\filldraw[black] (11, 0) coordinate (A2) circle (4pt) node{};
\filldraw[black] (11, 1) coordinate (A3) circle (4pt) node{};
\filldraw[black] (12, -1) coordinate (A4) circle (4pt) node{};
\filldraw[black] (12, 0) coordinate (A5) circle (4pt) node{};
\filldraw[black] (12, 2) coordinate (A6) circle (4pt) node{};
\filldraw[black] (13, -2) coordinate (A7) circle (4pt) node{};
\filldraw[black] (14, -3) coordinate (A8) circle (4pt) node{};    
\filldraw[black] (15, -4) coordinate (A9) circle (4pt) node{};
\filldraw[black] (9, 1) coordinate (A10) circle (4pt) node{};    
\filldraw[black] (8, 1) coordinate (A11) circle (4pt) node{};
 \filldraw[black] (14, -2) coordinate (A12) circle (4pt) node{};

%%%%% Draw edges %%%%%
\draw(A1) -- (A2) -- (A3) -- (A0) -- (A1);
\draw(A1) -- (A5);
\draw(A0) -- (A10) -- (A11);
\draw (A7) -- (A8) -- (A9);
\draw(A3) -- (A6);
\draw(A2)-- (A4) -- (A7) -- (A12);
\end{tikzpicture}
\caption{Two graphs with the same CSF on 13 vertices.}
\label{fig:n=13}
\end{figure}

\begin{figure}[hbt!]
\centering    
\begin{tikzpicture}[auto=center,every node/.style={circle, fill=black, scale=0.45}, style=thick, scale=0.4] 
%%%%% Draw vertices %%%%%
\filldraw[black] (0, 2) coordinate (A0) circle (4pt) node{};
\filldraw[black] (2, 0) coordinate (A1) circle (4pt) node{};
\filldraw[black] (1, -2) coordinate (A2) circle (4pt) node{};
\filldraw[black] (-1, -2) coordinate (A3) circle (4pt) node{};
\filldraw[black] (-2, 0) coordinate (A4) circle (4pt) node{};
\filldraw[black] (0, 3) coordinate (A5) circle (4pt) node{};
\filldraw[black] (3, -2) coordinate (A6) circle (4pt) node{};
\filldraw[black] (4, -2) coordinate (A7) circle (4pt) node{};
\filldraw[black] (5, -2) coordinate (A8) circle (4pt) node{};    
\filldraw[black] (6, -2) coordinate (A9) circle (4pt) node{};
\filldraw[black] (3, -1) coordinate (A10) circle (4pt) node{};    
\filldraw[black] (4, -1) coordinate (A11) circle (4pt) node{};

\filldraw[black] (-3, 0) coordinate (A12) circle (4pt) node{};
\filldraw[black] (-2, -2) coordinate (A13) circle (4pt) node{};
\filldraw[black] (-3, -2) coordinate (A14) circle (4pt) node{};
\filldraw[black] (-1, -3) coordinate (A15) circle (4pt) node{};
%%%%% Draw edges %%%%%
\draw(A5) -- (A0) -- (A1) -- (A2) -- (A3) -- (A4) -- (A0);
\draw(A2) -- (A6) -- (A7) -- (A8) -- (A9);
\draw(A6) -- (A10);
\draw(A7) -- (A11);
\draw(A4) -- (A12);
\draw(A15) -- (A3) -- (A13) -- (A14);
\end{tikzpicture}\qquad
\begin{tikzpicture}[auto=center,every node/.style={circle, fill=black, scale=0.45}, style=thick, scale=0.4] 
%%%%% Draw vertices %%%%%
\filldraw[black] (0, 2) coordinate (A0) circle (4pt) node{};
\filldraw[black] (2, 0) coordinate (A1) circle (4pt) node{};
\filldraw[black] (1, -2) coordinate (A2) circle (4pt) node{};
\filldraw[black] (-1, -2) coordinate (A3) circle (4pt) node{};
\filldraw[black] (-2, 0) coordinate (A4) circle (4pt) node{};
\filldraw[black] (0, 3) coordinate (A5) circle (4pt) node{};
\filldraw[black] (3, -2) coordinate (A6) circle (4pt) node{};
\filldraw[black] (4, -2) coordinate (A7) circle (4pt) node{};
\filldraw[black] (5, -2) coordinate (A8) circle (4pt) node{};    
\filldraw[black] (6, -2) coordinate (A9) circle (4pt) node{};
\filldraw[black] (3, -1) coordinate (A10) circle (4pt) node{};    
\filldraw[black] (4, -1) coordinate (A11) circle (4pt) node{};

\filldraw[black] (-3, 0) coordinate (A12) circle (4pt) node{};
\filldraw[black] (3, 0) coordinate (A13) circle (4pt) node{};
\filldraw[black] (-4, 0) coordinate (A14) circle (4pt) node{};
\filldraw[black] (1, -3) coordinate (A15) circle (4pt) node{};
%%%%% Draw edges %%%%%
\draw(A5) -- (A0) -- (A1) -- (A2) -- (A3) -- (A4) -- (A0);
\draw(A15) -- (A2) -- (A6) -- (A7) -- (A8) -- (A9);
\draw(A6) -- (A10);
\draw(A7) -- (A11);
\draw(A4) -- (A12) -- (A14);
\draw(A1) -- (A13);
\end{tikzpicture}
\caption{Two graphs with the same CSF on 16 vertices.}
\label{fig:n=16}
\end{figure}
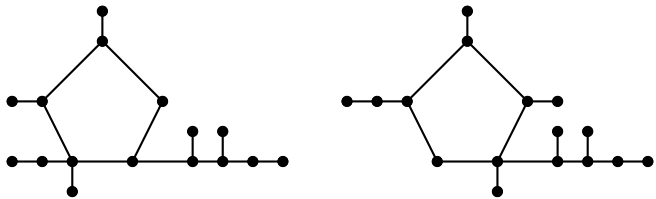

There is a large family of pairs of non-isomorphic connected unicyclic graphs with a 3-cycle and an even number of vertices that have identical CSFs but distinct numbers of non-trivial rooted trees, internal edges, and degree sequences. On the other hand, the following conjecture has been verified for $n\leq 16$ vertices.
\begin{conjecture}
    Let $G_1$ and $G_2$ be connected unicyclic graphs, either with a 3-cycle and an odd number of vertices or a $c$-cycle for $c\geq4$. Then if $\mathbf{X}_{G_1}=\mathbf{X}_{G_2}$ then $G_1$ and $G_2$ have the same degree sequence, number of non-trivial rooted trees, and number of internal edges.
\end{conjecture}
In \Cref{fig: n=12,fig:n=13,fig:n=16}, each pair of graphs shares the same degree sequence, number of non-trivial rooted trees, and number of internal edges. By \Cref{cor:cyc size hooks} two unicyclic graphs with the same CSF must have the same cycle size $c$. However, it is not known whether there exist two unicyclic graphs with different degree sequences that have the same CSF with $c > 3$.
We have verified the following using {\sc SageMath}~\cite{sagemath}.
\begin{enumerate}
\item CSF distinguishes $3$-unicyclic graphs on $n$ vertices if $n = 3,4,5,7,9,11,13,15$.
\item There are pairs of $3$-unicyclic graphs with the same CSF when the number of vertices is $ 6,8,10,12,14,16$.
\item CSF distinguishes $4$-unicyclic graphs on $n$ vertices if $n = 14, 15, 16, 17$.
\item CSF distinguishes $5$-unicyclic graphs on $n$ vertices if $n \leq 15$.
\item CSF distinguishes $c$-unicyclic graphs on $n$ vertices if $n \leq 16$ and $6 \leq c\leq 16$.
\end{enumerate}

\printbibliography

\end{document}